# OPTIMAL SCALINGS FOR LOCAL METROPOLIS–HASTINGS CHAINS ON NONPRODUCT TARGETS IN HIGH DIMENSIONS[1]

By Alexandros Beskos, Gareth Roberts and Andrew Stuart

*University of Warwick*

We investigate local MCMC algorithms, namely the random-walk Metropolis and the Langevin algorithms, and identify the optimal choice of the local step-size as a function of the dimension $n$ of the state space, asymptotically as $n \to \infty$. We consider target distributions defined as a change of measure from a product law. Such structures arise, for instance, in inverse problems or Bayesian contexts when a product prior is combined with the likelihood. We state analytical results on the asymptotic behavior of the algorithms under general conditions on the change of measure. Our theory is motivated by applications on conditioned diffusion processes and inverse problems related to the 2D Navier–Stokes equation.

**1. Introduction.** The Markov chain Monte Carlo (MCMC) methodology provides a flexible approach for simulating high-dimensional distributions appearing, for instance, as posterior information in a Bayesian context [15] or as integrators in importance sampling methods [14]. This paper is concerned with local MCMC algorithms, namely the random-walk Metropolis (RWM) and the Langevin algorithms. Our objective is to investigate the optimal choice of the local step-size as a function of the dimension $n$ of the state space, asymptotically as $n \to \infty$. In particular, we examine if the step-size should diminish with $n$ and, if so, at what rate.

The results in this paper extend significantly those in [21] and [22] since we here consider a family of *nonproduct* target laws. Our theory covers practical, involved probabilistic models: we will consider conditioned diffusion processes and inverse problems related to the 2D Navier–Stokes equation. In both these cases, the target measure is a change of measure from a Gaussian law on some appropriate infinite-dimensional Hilbert space. Gaussian

Received July 2008; revised September 2008.

[1]Supported by an EPSRC grant.

laws on abstract Hilbert spaces correspond to product Gaussian densities via the Karhunen–Loève expansion. This motivates the setting in which we present our results: target distributions defined as change of measure from a product law.

Consider the target probability density $\pi_n : \mathbb{R}^n \mapsto \mathbb{R}$ defined w.r.t. the Lebesque measure $dx$. The Metropolis–Hastings theory defines a Markov chain reversible w.r.t. $\pi_n$. Under the assumption of ergodicity, a path of the chain will provide, once convergence to stationarity is achieved, correlated draws from $\pi_n$. The dynamics of the chain are developed as follows. Given the current location $x \in \mathbb{R}^n$ a move to $y$ is proposed according to some user-specified transition kernel $q_n(x, dy) = q_n(x, y) \, dy$. Reversibility w.r.t. $\pi_n$ is then guaranteed [18] if the chain moves from $x \to x'$ according to the rule:

$$(1.1) \qquad x' = \begin{cases} y, & \text{with probability } a_n(x, y), \\ x, & \text{otherwise,} \end{cases}$$

where the acceptance probability is given by

$$(1.2) \qquad a_n(x, y) = \begin{cases} 1 \wedge \dfrac{\pi_n(y) q_n(y, x)}{\pi_n(x) q_n(x, y)}, & \text{if } \pi_n(x) q_n(x, y) > 0, \\ 0, & \text{if } \pi_n(x) q_n(x, y) = 0. \end{cases}$$

Ergodicity holds under regularity conditions on $\pi_n$, $q_n$; see, for instance, [15].

In this paper, $\pi_n$ will be a change of measure from a reference product law, denoted by $\tilde{\pi}_n$. In particular, the family of targets is specified as follows:

$$(1.3a) \qquad \frac{d\pi_n}{d\tilde{\pi}_n}(x) = \exp\{-\phi_n(x)\}$$

for some $\mathcal{B}(\mathbb{R}^n)$-measurable mapping $\phi_n$, with $\tilde{\pi}_n$ having Lebesque density:

$$(1.3b) \qquad \tilde{\pi}_n(x) = \prod_{i=1}^{n} \frac{1}{\lambda_i} f\left(\frac{x_i}{\lambda_i}\right) = \prod_{i=1}^{n} \frac{1}{\lambda_i} \exp\{-g(x_i/\lambda_i)\}$$

for appropriate $f, g : \mathbb{R} \mapsto \mathbb{R}$. Motivated by applications, the standard deviations are assumed to satisfy

$$(1.3c) \qquad \lambda_i = i^{-\kappa}, \qquad i = 1, 2, \ldots,$$

for some $\kappa \geq 0$.[2] We investigate MCMC algorithms corresponding to the following local-move proposals:

$$(1.4) \qquad \text{RWM: } y = x + \sigma_n Z,$$

$$(1.5) \qquad \text{SLA: } y = x + \frac{\sigma_n^2}{2} \nabla \log \tilde{\pi}_n(x) + \sigma_n Z$$

---

[2]We later relax this assumption to allow for algebraic decay at rate $i^{-\kappa}$; see (2.6).

with $Z \sim \mathcal{N}(0, I_n)$ and step-size $\sigma_n > 0$. The first corresponds to a random-walk update; the second to a simplified version of the standard Metropolis-adjusted Langevin algorithm (MALA) with proposal:

$$(1.6) \qquad \text{MALA:} \ y = x + \frac{\sigma_n^2}{2} \nabla \log \pi_n(x) + \sigma_n Z.$$

We will show that using only the reference law $\tilde{\pi}_n$ in SLA (which stands for "simplified Langevin algorithm") does not reduce the asymptotic properties of the resulting algorithm. The objective is the identification of the "appropriate" scaling of $\sigma_n$ as a function of the dimension $n$, for large $n$. Excessively large or small $\sigma_n$ will typically result in unacceptably small or high acceptance probabilities, respectively, and to poor ergodicity and mixing properties for the MCMC.

References [21] and [22] analyze RWM and MALA as applied to i.i.d. targets corresponding to our $\tilde{\pi}_n$ with $\kappa = 0$. They show that step-sizes should be arranged as $\sigma_n^2 = \mathcal{O}(n^{-1})$ for RWM and $\sigma_n^2 = \mathcal{O}(n^{-1/3})$ for MALA. Extensions of such results have thereafter appeared in several papers; see [1, 2] and [20] for product targets. For nonproduct targets, [8] examines an exchangeable target, [9] a local interaction model related with Gibbs distributions and [24] elliptically symmetric targets. See [3] for an analytical review. Using a different approach, we look at the structure of the nonproduct target $\pi_n$ in (1.3) and present general, analytical conditions on $\phi_n$ under which we show that one should select $\sigma_n^2 = \mathcal{O}(n^{-2\kappa-1})$ for RWM and $\sigma_n^2 = \mathcal{O}(n^{-2\kappa-1/3})$ for SLA. We will use the average squared-jump-distance as an index of the efficiency of MCMC algorithms, as it allows for transparent, explicit calculations. Our analysis is considerably simpler than the approach adopted in the preceding papers and, as we will show, the results are relevant for probability measures arising in practical applications.

A motivation for investigating the change of measure (1.3) are cases when the densities $\pi_n$, $\tilde{\pi}_n$ under consideration are finite-dimensional approximations of some *infinite*-dimensional measures $\pi$, $\tilde{\pi}$ related through the density

$$(1.7) \qquad \frac{d\pi}{d\tilde{\pi}}(X) = \exp\{-\phi(X)\}$$

on an appropriate space. In this paper we substantiate the intuition that the presence of absolute continuity in the limit $n \to \infty$ implies similar asymptotic behavior for the MCMC between the product and nonproduct scenarios.

The structure of the paper is as follows. In Section 2 we state and prove the scaling results for RWM and SLA. In Sections 3, 4 and 5 we present applications of the theory in the context of conditioned diffusion processes and Navier–Stokes inverse problems. In Section 6 we show some additional ideas and algorithms which can deliver improved results. We finish with some conclusions in Section 7. Proofs are collected in an Appendix B; before that, we have collected some Taylor expansion results needed in the proofs in Appendix A.

**2. The main results.** We will state rigorously the scaling results outlined in the Introduction. For both RWM and SLA, as applied to our target $\pi_n$ in (1.3), the acceptance probability in (1.2) can be written as

$$(2.1) \qquad a_n(x, y) = 1 \wedge \exp\{\phi_n(x) - \phi_n(y) + R_n(x, y)\}$$

for an exponent $R_n = R_n(x, y)$, which in the case of RWM is equal to

$$(2.2) \qquad R_n^{\mathrm{RWM}} = \sum_{i=1}^{n} \{g(x_i/\lambda_i) - g(y_i/\lambda_i)\},$$

whereas for SLA:

$$
\begin{aligned}
(2.3) \qquad R_n^{\mathrm{SLA}} &= R_n^{\mathrm{RWM}} + \sum_{i=1}^{n} \delta_{i,n}, \\
\delta_{i,n} &:= \left(\frac{y_i}{\lambda_i} - \frac{x_i}{\lambda_i}\right) \frac{g'(y_i/\lambda_i) + g'(x_i/\lambda_i)}{2} \\
&\qquad - \frac{\sigma_n^2}{8\lambda_i^2}((g'(y_i/\lambda_i))^2 - (g'(x_i/\lambda_i))^2).
\end{aligned}
$$

In the case of $R_n^{\mathrm{RWM}}$, $y$ and $x$ are connected via (1.4), in the case of $R_n^{\mathrm{SLA}}$ via (1.5). We impose the following regularity conditions on $g$ and the density $f = \exp\{-g\}$.

CONDITION 1. (i) The function $g$ is infinitely differentiable with derivatives of all orders having a polynomial growth bound.

(ii) All moments of $f$ are finite.

We henceforth assume that Condition 1 holds without further reference. We will explore the behavior of MCMC algorithms in equilibrium; we note that algorithmic properties could be different in transient regimes; see [10]. Thus, most expectations are considered in stationarity and $E[\cdot]$ will in general denote expectation under the relevant target $\pi_n$; we will sometimes write $E_{\pi_n}[\cdot]$ or $E_{\tilde{\pi}_n}[\cdot]$ when we need to be explicit about the measure under consideration.

Broadly speaking, our approach will be to identify the limiting properties of the exponent $\phi_n(x) - \phi_n(y) + R_n$ as $n \to \infty$. Through this term we will then be able to derive the asymptotic properties of the average acceptance probability $E_{\pi_n}[a_n(x, y)]$ and obtain analytical results for the squared-jump-distance over the various choices of the step-size $\sigma_n$. We will work with $L_q$-limits and convergence in law, and we will obtain results under fairly general assumptions, suited to complex (nonproduct) probabilistic models arising in applications. In contrast, the analysis so far in the literature has been somewhat technical, focusing most on a rich array of results available in the (restrictive) product set-up.

2.1. *The product case.* To connect with the literature, we will initially state a result for the simplified product case when $\phi_n \equiv 0$ for all $n \geq 1$. Several of the ideas in the proof will be relevant in the nonproduct scenario later. All stated expectations are in stationarity (here $x \sim \tilde{\pi}_n$); $\Phi$ is the CDF of $\mathcal{N}(0,1)$. To avoid repetition when stating results for RWM and SLA we define the following index $I$:

$$\text{RWM: } I = 1,$$

(2.4)

$$\text{SLA: } I = 1/3.$$

THEOREM 1. *Let $\tilde{\pi}_n$ in (1.3) be the target distribution for RWM and SLA. Assume that $\sigma_n^2 = l^2 n^{-\rho}$ for $l, \rho > 0$. Then, as $n \to \infty$:*

(i) *if $\rho = 2\kappa + I$, for $I$ as in (2.4), then $E[a_n(x,y)] \to a(l)$, where*

$$a^{\mathrm{RWM}}(l) = 2\Phi\left(-\frac{l}{2}\sqrt{\frac{K^{\mathrm{RWM}}}{\tau^{\mathrm{RWM}}}}\right), \qquad a^{\mathrm{SLA}}(l) = 2\Phi\left(-\frac{l^3}{2}\sqrt{\frac{K^{\mathrm{RWM}}}{\tau^{\mathrm{SLA}}}}\right),$$

*for constants $K^{\mathrm{RWM}}, K^{\mathrm{SLA}} > 0$, with $\tau^{\mathrm{RWM}} = 1 + 2\kappa$ and $\tau^{\mathrm{SLA}} = 1 + 6\kappa$,*

(ii) *if $\rho > 2\kappa + I$, then $E[a_n(x,y)] \to 1$,*

(iii) *if $2\kappa < \rho < 2\kappa + I$, then $n^p E[a_n(x,y)] \to 0$ for any $p \geq 0$.*

The constants $K^{\mathrm{RWM}}$ and $K^{\mathrm{SLA}}$ are given in the proof and depend only on the density $f$ appearing at the definition of $\tilde{\pi}_n$. The results are based on the limiting behavior of $R_n$. When $\rho = 2\kappa + I$, $R_n$ is subject to a central limit theorem, forcing a limit also for the average acceptance probability. Other step-sizes will $n$-eventually lead to a degenerate acceptance probability. We will see in the sequel that a quantity providing an index for the efficiency of the MCMC algorithms considered here, is the product

$$\sigma_n^2 E[a_n(x,y)].$$

Clearly, in the context of Theorem 1, this quantity is maximized (in terms of the dimension $n$) for $\rho = 2\kappa + I$; it is larger for SLA than for RWM due to the presence of information about the target in the proposal.

2.2. *The general case.* The above results describe the behavior of the expectation $E_{\tilde{\pi}_n}[1 \wedge \exp\{R_n\}]$ for large $n$. One could now look for conditions on the change of measure (1.3a) that allow some of these results to apply to the more interesting nonproduct scenario. The quantity under investigation is now $E_{\pi_n}[1 \wedge \exp\{\phi_n(x) - \phi_n(y) + R_n\}]$, an expectation w.r.t. $\pi_n$. We will first consider the following condition:

CONDITION 2.   There exists $M > 0$ such that for any sufficiently large $n$

$$|\phi_n(x)| \leq M \qquad \text{for all } x \in \mathbb{R}^n.$$

Of course, it would suffice that the condition holds $\tilde{\pi}_n$-a.s.; in the sequel, conditions on $\phi_n$ should be interpreted as conditions on some $\tilde{\pi}_n$-version of $\phi_n = -\log(d\pi_n/d\tilde{\pi}_n)$.

THEOREM 2.   Let $\pi_n$ in (1.3) be the target distribution for RWM and SLA. Assume that $\sigma_n^2 = l^2 n^{-\rho}$ for $l, \rho > 0$. If $\{\phi_n\}$ satisfies Condition 2, then as $n \to \infty$:

   (i) if $\rho \geq 2\kappa + I$, then $\liminf_{n\to\infty} E[a_n(x,y)] > 0$,
   (ii) if $2\kappa < \rho < 2\kappa + I$, then $n^p E[a_n(x,y)] \to 0$ for any $p \geq 0$,

for the appropriate index $I$ for each of RWM and SLA specified in (2.4).

Condition 2 provides a recipe for a direct transfer of some of the results of the product case to the nonproduct one. A more involved result, motivated by the collection of applications we describe in the following sections, is given in the next theorem. For any $s \in \mathbb{R}$, we define the norm $|\cdot|_s$ on $\mathbb{R}^n$ as follows:

$$(2.5) \qquad |x|_s = \left(\sum_{i=1}^n i^{2s} x_i^2\right)^{1/2}.$$

We also set $|\cdot| \equiv |\cdot|_0$.

CONDITION 3.   There exists $M \in \mathbb{R}$ such that, for all sufficiently large $n$,

$$\phi_n(x) \geq M \qquad \text{for all } x \in \mathbb{R}^n.$$

CONDITION 4.   There exist constants $C > 0$, $p > 0$ and $s < \kappa - 1/2$, $s' < \kappa - 1/2$ such that, for all sufficiently large $n$,

$$|\phi_n(y) - \phi_n(x)| \leq C(1 + |x|_s^p + |y - x|_s^p)|y - x|_{s'}$$

for all $x, y \in \mathbb{R}^n$.

Condition 4 is motivated by the application to conditioned diffusions in one of the following sections; we shall see that, in some contexts, it follows from a polynomial growth assumption on the derivative of $\phi_n$. Condition 3 prevents the nonproduct measure $\pi_n$ from charging sets with much higher probability than $\tilde{\pi}_n$. The following condition (clearly weaker than Condition 4) will be relevant in the Navier–Stokes problem. We set $\mathbb{R}^+ = [0, \infty)$.

CONDITION 5. There exist locally bounded function $\delta : \mathbb{R}^+ \times \mathbb{R}^+ \mapsto \mathbb{R}^+$ and constants $C > 0$, $p > 0$, and $s, s', s''$ all three (strictly) smaller than $\kappa - 1/2$, such that, for all sufficiently large $n$:

(i) $|\phi_n(y) - \phi_n(x)| \leq \delta(|x|_s, |y|_s)|y - x|_{s'}$,
(ii) $|\phi_n(x)| \leq C(1 + |x|_{s''}^p)$ for all $x, y \in \mathbb{R}^n$.

THEOREM 3. Let $\pi_n$ in (1.3) be the target distribution for RWM and SLA. Assume that $\sigma_n^2 = l^2 n^{-\rho}$ for $l, \rho > 0$. If $\{\phi_n\}$ satisfies Conditions 3 and 5, then as $n \to \infty$:

(i) if $\rho = 2\kappa + I$, then $E[a_n(x,y)] \to a(l)$, for $a(l)$ defined in Theorem 1,
(ii) if $\rho > 2\kappa + I$, then $E[a_n(x,y)] \to 1$,
(iii) if $2\kappa < \rho < 2\kappa + I$, then $n^p E[a_n(x,y)] \to 0$ for any $p \geq 0$.

REMARK 1. Condition 5 implies the probabilistic statement

$$\phi_n(y) - \phi_n(x) \xrightarrow{L_q(\tilde{\pi}_n)} 0$$

for arbitrarily large $q$. This result, together with Condition 3, suffices for showing that the effect of the change of measure in the scaling properties of RWM and SLA is asymptotically negligible.

2.3. *Optimality.* We use the squared-jump-distance as an index for the efficiency of different MCMC algorithms. Specifically, for the algorithms we have considered so far, we will calculate the quantity

$$S_n := E[(x'_{i^*} - x_{i^*})^2]$$

for some arbitrary fixed $i^*$; here $x'$ is the location of the MCMC Markov chain after one step given that currently the chain is at $x \sim \pi_n$ and $i^*$ refers to a fixed element of $\{1, 2, \ldots, n\}$. Note that

$$\text{Corr}_n(x_{i^*}, x'_{i^*}) = 1 - \frac{S_n}{2\text{Var}_n},$$

$\text{Corr}_n$ denoting correlation and $\text{Var}_n$ the variance of $x_{i^*}$ in stationarity. Thus, larger $S_n$ implies lower first-order autocorrelation.

THEOREM 4. Let $\pi_n$ in (1.3) be the target distribution for RWM and SLA. Assume that $\sigma_n^2 = l^2 n^{-\rho}$ for $l, \rho > 0$. If $\{\phi_n\}$ satisfies Conditions 3 and 5, then as $n \to \infty$:

(i) if $\rho = 2\kappa + I$, then

$$S_n = l^2 a(l) \times n^{-2\kappa - I} + o(n^{-2\kappa - I})$$

for $a(l)$ defined in Theorem 1,

(ii) *if $\rho > 2\kappa + I$, then $S_n = l^2 n^{-\rho} + o(n^{-\rho})$,*
(iii) *if $2\kappa < \rho < 2\kappa + I$, then $S_n = \mathcal{O}(n^{-p})$ for any $p \geq 0$.*

Theorem 4 shows how to scale the proposal step, as a function of $n$, to maximize the mean squared-jump-distance; we can then tune the parameter $l$. When maximizing the coefficient $l^2 a(l)$ over $l$ we retrieve a familiar result characterizing RWM and the Langevin algorithms: in the case of RWM, $l^2 a(l)$ is maximized for that $l$ for which $a^{\text{RWM}}(l) = 0.234$, and in the case of SLA for $a^{\text{SLA}}(l) = 0.574$, for *any* choice of the reference density $f$ and any change of measure satisfying Conditions 3 and 5. These characteristic numbers were first obtained in [21] and [22], in the simplified context of i.i.d. target distributions.

2.4. *A generalization.* It is straightforward to extend the results so far stated to a larger family of target distributions. Such a generalization will be required in the applications considered in the sequel.

COROLLARY 1. *Allow the target distribution $\pi_n$ to be as in (1.3) but with parameters $\lambda_i$ that may also depend on $n$, so $\lambda_i = \lambda_{i,n}$. Assume that there exist constants $C_-, C_+ > 0$ and $\kappa \geq 0$ such that, for all $n$ and $1 \leq i \leq n$,*

$$(2.6) \qquad C_- i^{-\kappa} \leq \lambda_{i,n} \leq C_+ i^{-\kappa}.$$

*Then all statements of Theorems 1–4 apply directly to the target $\pi_n$ with $\lambda_i = \lambda_{i,n}$ except for those where the limiting probabilities $a^{\text{RWM}}(l)$ and $a^{\text{SLA}}(l)$ appear. For the latter cases, the statements will hold after replacing the stated values for $\tau^{\text{RWM}}$ and $\tau^{\text{SLA}}$ with*

$$\tau^{\text{RWM}} = \lim_{n \to \infty} n^{-(2\kappa+1)} \sum_{i=1}^{n} \lambda_{i,n}^{-2}, \qquad \tau^{\text{SLA}} = \lim_{n \to \infty} n^{-(6\kappa+1)} \sum_{i=1}^{n} \lambda_{i,n}^{-6},$$

*provided that the above limits exist.*

**3. Applications in infinite dimensions.** In this section we introduce a range of applications which require the sampling of a probability measure on function space. The common mathematical structure of these problems is that the measure of interest, $\pi$, has density with respect to a Gaussian reference measure $\tilde{\pi} \sim \mathcal{N}(m, \mathcal{C})$:

$$(3.1) \qquad \frac{d\pi}{d\tilde{\pi}}(X) \propto \exp\{-\phi(X)\}.$$

The Karhunen–Loève expansion for Gaussian measures on a Hilbert space allows us to view the target $\pi$ as a change of measure from a product of scalar Gaussian densities on $\mathbb{R}^\infty$, thus casting the simulation problem into the theory presented in the first part of the paper.

We highlight two types of problem where the structure (3.1) arises naturally. The first type concerns SDEs conditioned on observations. Here the Gaussian reference measure is typically a conditioned Gaussian diffusion, in which nonlinear drifts are ignored, and the posterior measure is found via application of the Girsanov formula. The second type of problem concerns inverse problems for differential equations, where prior knowledge about an unknown function, in the form of a Gaussian measure, is combined with observations, via application of the Bayes formula, to determine a posterior measure on function space. The common structure inherent in these problems allows for the use of the same notation in the different contexts: the mean of the reference Gaussian measure will be $m$ and its covariance operator will be $\mathcal{C}$; the *precision operator* $-\mathcal{C}^{-1}$ will be denoted by $\mathcal{L}$. The state space will be a separable Hilbert space $\mathcal{H}$.

The Gaussian law is well defined if and only if $\mathcal{C}: \mathcal{H} \mapsto \mathcal{H}$ is a positive, self-adjoint and trace-class operator, the last property meaning that its eigenvalues are summable. Thus, we may construct an orthonormal basis of eigenfunctions $\{e_i\}_{i=1}^{\infty}$ and corresponding eigenvalues $\{\lambda_i^2\}_{i=1}^{\infty}$ satisfying $\mathcal{C}e_i = \lambda_i^2 e_i$ with $\sum_i \lambda_i^2 < \infty$. A typical draw $X \sim \mathcal{N}(m, \mathcal{C})$ can be written via the Karhunen–Loève expansion as

$$(3.2) \qquad X = \sum_{i=1}^{\infty} x_i e_i, \qquad x_i = m_i + \lambda_i \xi_i,$$

where the sequence $\{\xi_i\}_{i=1}^{\infty}$ is i.i.d. with $\xi_1 \sim \mathcal{N}(0,1)$ and $m_i = \langle m, e_i \rangle$. It is readily verified that

$$E[(X - m) \otimes (X - m)] = \sum_{i=1}^{\infty} \lambda_i^2 (e_i \otimes e_i) = \mathcal{C},$$

which agrees (here $\otimes$ stands for tensor product, and expectation is over a random linear operator; see [12] for more details on theory of expectations on general Hilbert spaces) with the familiar identity for the covariance matrix on Euclidean spaces. From (3.2), the isomorphism $X \mapsto \{x_i\}$ allows us to view $\mathcal{N}(m, \mathcal{C})$ as a product measure on $\ell_2$, the space of square summable sequences.

A Gaussian measure on a Hilbert space is thus easy to build, given an orthonormal basis for the space, simply by specifying the eigenvalues of the covariance operator. Such an approach provides also a natural way to specify regularity properties of Hilbert space valued random functions. To illustrate this, we define the norm $\| \cdot \|_s$ on $\mathcal{H}$ (which by the duality $X \mapsto \{x_i\}$ can also be seen as a norm in $\ell_2$) as follows:

$$(3.3) \qquad \|X\|_s = \left( \sum_{i=1}^{\infty} i^{2s} x_i^2 \right)^{1/2}$$

for $s \in \mathbb{R}$. The largest $s$ for which $\|X\|_s$ is finite is a measure of the regularity of $X$ as it encodes information about the rate of decay of the coefficients $x_i$. From the Karhunen–Loève expansion we deduce that

$$E\|X - m\|_s^2 = \sum_{i=1}^{\infty} i^{2s} \lambda_i^2.$$

Thus, the Gaussian measure delivers realizations of finite $\|\cdot\|_s$-norm if $\lambda_i^2$ is chosen to decay like $i^{-(2s+1+\varepsilon)}$ for some $\varepsilon > 0$. Depending on the particular space under consideration, the definition of the $\|\cdot\|_s$ can slightly vary from the one in (3.3); we will give explicit definitions for each of the applications in the sequel. In many cases there is a natural relationship between the spaces of finite $\|\cdot\|_s$-norm and various function spaces which arise naturally in the theory of partial differential equations. For example, in one dimension with function $e_i(t) = \sqrt{2}\sin(i\pi t)$ one obtains Sobolev spaces (see, e.g., [23]).

REMARK 2. Finite-dimensional densities $\pi_n$, $\tilde{\pi}_n$, corresponding to approximations of $\pi$, $\tilde{\pi}$, can be derived by spectral methods, finite difference methods or other approaches. We present such directions in the context of the example applications that follow. We show that under general conditions on the concrete functional $\phi$ under investigation, its discretized counterpart $\phi_n$ will satisfy some of the conditions introduced in Section 2. Thus, we can use the results from that section to optimize MCMC methods applied to the complicated sampling problems presented in the next two sections.

**4. Diffusion bridges.** There is a variety of applications where it is of interest to study SDEs conditioned to connect two points in space-time: the so-called *diffusion bridges*. We limit our study to problems with additive noise satisfying the equation

$$\frac{dX}{dt} = h(X) + \gamma \frac{dW}{dt},$$

$$X(0) = x^-, \qquad X(T) = x^+.$$

Here $h : \mathbb{R}^d \mapsto \mathbb{R}^d$ and $\gamma \in \mathbb{R}^{d \times d}$; we define $\Gamma = \gamma\gamma^T$ and assume that $\Gamma$ is invertible. We denote by $\tilde{\pi}$ the $d$-dimensional Brownian bridge measure arising when $h \equiv 0$ and by $\pi$ the general non-Gaussian measure. Under mild conditions on $h$ (see, e.g., [13]), the two measures are absolutely continuous w.r.t. each other. Their Radon–Nikodym derivative is given by Girsanov's theorem and assumes the general form (3.1) where now

$$\phi(X) = \int_0^T \left( \frac{1}{2}|h(X)|_\Gamma^2 \, dt - \langle h(X), dX \rangle_\Gamma \right)$$

with $\langle \cdot, \cdot \rangle_\Gamma = \langle \cdot, \Gamma^{-1} \cdot \rangle$, the latter being the standard inner product on $\mathbb{R}^d$. In this context $\mathcal{H} = L_2([0, T], \mathbb{R}^d)$. The expression for $\phi$ above involves a stochastic integral; it can be replaced with a Riemann integral in the particular scenario when $\Gamma^{-1} h$ is a gradient. In this context we will consider the bridge diffusion:

$$\frac{dX}{dt} = -\nabla V(X) + \sqrt{\frac{2}{\beta}} \frac{dW}{dt},$$

(4.1)

$$X(0) = x^-, \qquad X(T) = x^+$$

with inverse temperature $\beta > 0$ and potential $V : \mathbb{R}^d \mapsto \mathbb{R}$. Applying Itô's formula and ignoring constants, we get that

$$\phi(X) = \int_0^T G(X) \, dt,$$

(4.2)

where

$$G(u) = \frac{\beta}{4} |\nabla V(u)|^2 - \frac{1}{2} \Delta V(u), \qquad u \in \mathbb{R}^d.$$

A typical application from molecular dynamics is illustrated in Figure 1. The figure shows a crystal lattice of atoms in two dimensions, with an atom removed from one site. The potential $V$ here is a sum of pairwise potentials between atoms which has an $r^{-12}$ repulsive singularity. The lattice should be viewed as spatially extended to the whole of $\mathbb{Z}^2$ by periodicity. The removal of an atom creates a vacancy which, under thermal activation as in (4.1), will diffuse around the lattice: the vacancy moves lattice sites whenever one of the neighboring atoms moves into the current vacancy position. This motion of the atoms is a rare event, and we can condition our model on its occurrence. This application, as well as others involving conditioned diffusions arising in signal processing, are detailed in [6].

4.1. *The Fourier expansion.* We focus on model (4.1). The mean of the reference Brownian bridge Gaussian measure is $m(t) = x^-(1 - t/T) + x^+ t/T$. Without loss of generality we will assume that $x^- = x^+ = 0$, therefore $m \equiv 0$; otherwise, one should work with $X - m$. The precision operator of the Brownian bridge is the Laplacian (see [16]):

$$\mathcal{L} = \frac{\beta}{2} \frac{d^2}{dt^2}$$

with Dirichlet boundary conditions on $[0, T]$, specified through the domain of $\mathcal{L}$. In the scalar case $d = 1$, the (orthonormal) eigenfunctions $\{e_i\}$ and eigenvalues $\{\lambda_i^2\}$ of the covariance operator $\mathcal{C} = -\mathcal{L}^{-1}$ are

$$e_i = e_i(t) = \sqrt{\frac{2}{T}} \sin(i\pi t/T), \qquad \lambda_i^2 = \frac{2}{\beta} \frac{T^2}{\pi^2} i^{-2}$$

(4.3)

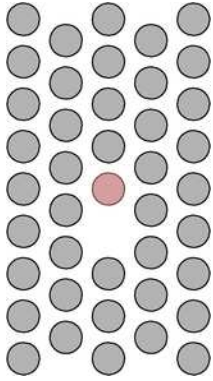

FIG. 1. *Crystal lattice with vacancy. We condition on the red atom moving into the vacancy.*

with $i \geq 1$. For general dimension $d$, the covariance operator is simply the diagonal operator matrix with copies of the covariance operator for $d = 1$ in the diagonal. In this context we will use the expansion

$$(4.4) \qquad X = \sum_{i=1}^{\infty} x_{\cdot,i} e_i, \qquad x_{\cdot,i} = (x_{1,i}, x_{2,i}, \ldots, x_{d,i})^{\top}.$$

A finite-dimensional density approximating $\exp\{-\phi(X)\}$ can be obtained via spectral truncation. For $n = d \times N$ and the vector in $\mathbb{R}^n$

$$(4.5) \qquad x = (x_{\cdot,1}^{\top}, x_{\cdot,2}^{\top}, \ldots, x_{\cdot,N}^{\top}),$$

we get the density[3]

$$(4.6a) \qquad \frac{d\pi_n}{d\tilde{\pi}_n}(x) \propto \exp\{-\phi(P_N(x))\}, \qquad P_N(x) := \sum_{i=1}^{N} x_{\cdot,i} e_i.$$

By the Karhunen–Loève expansion, the reference measure will be

$$(4.6b) \qquad \tilde{\pi}_n = \prod_{i=1}^{N} \left\{ \prod_{j=1}^{d} \mathcal{N}(0, \lambda_i^2) \right\}$$

for the $\lambda_i$'s in (4.3).

4.2. *The results.* For $s \in \mathbb{R}$, we define the norm:

$$\|X\|_s = \left( \sum_{i=1}^{\infty} i^{2s} |x_{\cdot,i}|^2 \right)^{1/2}.$$

We need the following condition:

---

[3] Note that $P_N(x)$ is the $L_2$-projection of $X$ onto the linear span of $\{e_1, e_2, \ldots, e_N\}$.

CONDITION 6. $G$ is bounded from below and there exist constants $C > 0$, $p > 0$ for which

$$|G(w) - G(v)| \leq C(1 + |v|^p + |w - v|^p)|w - v| \qquad \text{for all } v, w \in \mathbb{R}^d.$$

Note that the stated condition on $|G(w) - G(v)|$ is implied by a polynomial growth assumption on the gradient of $G$. Condition 6 will imply Conditions 3 and 4 making possible the application of the results of the previous section for the target $\pi_n$ in (4.6). To see this, for $X, Y \in \mathcal{H}$ we use Condition 6 to get

$$
\begin{aligned}
(4.7) \qquad |\phi(Y) - \phi(X)| &\leq C \int_0^T (1 + |X(t)|^p + |Y(t) - X(t)|^p)|Y(t) - X(t)| \, dt \\
&\leq C'(1 + \|X\|_{L_{2p}}^p + \|Y - X\|_{L_{2p}}^p)\|Y - X\|_{L_2},
\end{aligned}
$$

where for the second calculation we used Cauchy–Schwarz and the triangle inequality for the $L_2$-norm. Now, the Sobolev embedding theorem (stated in Lemma B.3 in Appendix B) gives that for any $2 \leq q < \infty$ one can select $s = s(q) < 1/2$ so that

$$\|X\|_{L_q} \leq C\|X\|_s$$

for all $X$. Thus, continuing from (4.7):

$$(4.8) \qquad |\phi(Y) - \phi(X)| \leq C(1 + \|X\|_s^p + \|Y - X\|_s^p)\|Y - X\|_{L_2}$$

for some $s = s(p) < 1/2$. Equation (4.8) can now be used to deduce that $\phi_n(x) = \phi(P_N(x))$ satisfies Condition 4.

PROPOSITION 1. *If $G$ satisfies Condition 6, then the limiting results of Theorems 3 and 4 apply for the case of $\pi_n$ in (4.6) under the specifications:*

$$\kappa = 1, \qquad \tau^{\mathrm{RWM}} = \frac{\beta\pi^2}{6T^2 d^2}, \qquad \tau^{\mathrm{SLA}} = \frac{\beta^3\pi^6}{56T^6 d^6}.$$

4.3. *Finite differences.* Similar results also hold for other approximation schemes such as finite differences. In particular, assuming that paths are approximately constant on subintervals of length $\Delta t = \frac{T}{N+1}$, we get that

$$\frac{d\pi_n}{d\tilde{\pi}_n}(X) \propto \exp\left\{-\sum_{t=1}^N G(X_t)\Delta t\right\}, \qquad \Delta t = \frac{T}{N+1}$$

for argument $X = (X_t)_{t=1}^N$, with $X_t \in \mathbb{R}^d$. The reference law $\tilde{\pi}_n$ represents the finite-dimensional distributions of $d$ independent Brownian bridges at the discrete-time instances $\Delta t, 2\Delta t, \ldots, N\Delta t$:

$$\tilde{\pi}_n = \prod_{i=1}^{d} \mathcal{N}(0, -L^{-1}\Delta t^{-1}), \qquad L = \frac{\beta}{2} \frac{1}{\Delta t^2} \begin{pmatrix} -2 & 1 & & & \\ 1 & -2 & 1 & & \\ & & \ddots & & \\ & & 1 & -2 & 1 \\ & & & 1 & -2 \end{pmatrix}.$$

The eigenvectors $\{e_i\}_{i=1}^{N}$ of $C = -L^{-1}$ and its eigenvalues $\{\lambda_{i,n}^2\}_{i=1}^{N}$ are as follows [in this context $e_i = (e_{i,1}, \ldots, e_{i,N})$]:

$$e_{i,t} = \sqrt{\frac{2}{T}} \sin\left(\frac{i\pi t}{N+1}\right), \qquad \lambda_{i,n}^2 = \frac{T^2}{2\beta} \left(\sin\left(\frac{i\pi}{2(N+1)}\right)(N+1)\right)^{-2}.$$

The natural inner product for $\mathbb{R}^N$ in this context is $\langle u, v \rangle = \sum_{i=1}^{N} u_i v_i \Delta t$ under which one can verify that $\{e_i\}_{i=1}^{N}$ is an orthonormal basis. We factorize $\tilde{\pi}_n$ via the Karhunen–Loève expansion. That is, we write

$$(4.9) \qquad X_t = \sum_{i=1}^{N} x_{\cdot,i} e_{i,t}$$

and work on the space of $x = (x_{\cdot,1}^\top, x_{\cdot,2}^\top, \ldots, x_{\cdot,N}^\top)$ when the target can be rewritten as

$$(4.10) \qquad \frac{d\pi_n}{d\tilde{\pi}_n}(x) \propto \exp\left\{ -\sum_{t=1}^{N} G(X_t)\Delta t \right\}, \qquad \tilde{\pi}_n = \prod_{i=1}^{N}\left\{ \prod_{j=1}^{d} \mathcal{N}(0, \lambda_{i,n}^2) \right\}.$$

PROPOSITION 2. *If $G$ satisfies Condition 6, then the limiting results of Theorems 3 and 4 apply for the case of $\pi_n$ in (4.10) under the specifications:*

$$\kappa = 1, \qquad \tau^{\mathrm{RWM}} = \frac{\beta}{T^2 d^2}, \qquad \tau^{\mathrm{SLA}} = \frac{5\beta^3}{2T^6 d^6}.$$

**5. Data assimilation.** Data assimilation is concerned with the optimal blending of observational data and mathematical model to enhance predictive capability. It has had major impact in the realm of weather forecasting and is increasingly used by oceanographers. As a concrete (and simplified) model of these applications from the area of fluid mechanics, we consider the problem of determining the initial condition for Navier–Stokes equations for the motion of an incompressible Newtonian fluid in a two-dimensional unit box, with periodic boundary conditions, given observations. The equation, determining the velocity $X = X(t, v)$ over the torus domain $\mathbb{T}^2$, can be written as follows:

$$\frac{dX}{dt} = \nu AX + B(X, X) + h, \qquad t > 0,$$

$$X(0, v) = X_0(v), \qquad\qquad\quad t = 0.$$

This should be rigorously interpreted as an ordinary differential equation in the Hilbert space $\mathcal{H}$, defined as the closure in $L^2(\mathbb{T}^2, \mathbb{R}^2)$ of the space of periodic, divergence-free and smooth functions on $[0, 1] \times [0, 1]$, with zero average (see, e.g., [11] for more details). We specify the operator $A$ in detail below, noting here that it is essentially the Laplacian on $\mathcal{H}$; the operator $B$ is a bilinear form and $h$ a forcing function, but details of these terms will not be relevant in what follows.

As a first example of the type of data encountered in weather forecasting we assume that we are given noisy observations of the velocity field at positions $\{v_l\}_{l=1}^L$ and times $\{t_m\}_{m=1}^M$. To be precise, we observe $\{W_{l,m}\}$ given by

$$W_{l,m} = X(v_l, t_m) + \xi_{l,m}, \qquad l = 1, \ldots, L, \ m = 1, \ldots, M,$$

where the $\xi_{l,m}$'s are zero-mean Gaussian random variables. Concatenating data we may write

$$W = \tilde{z} + \xi$$

for $W = (W_{1,1}, W_{2,1}, \ldots, W_{L,M})$, $\tilde{z} = (X(v_1, t_1), X(v_2, t_1), \ldots, X(v_L, t_M))$ and $\xi \sim \mathcal{N}(0, \Sigma)$ for some covariance matrix $\Sigma$ capturing the correlations in the noise. We refer to this setting as *Eulerian data assimilation*, the word Eulerian denoting the fact that the observations are of the Eulerian velocity field.

For a second example, illustrating a data type commonly occurring in oceanography, we assume that we are given noisy observations of Lagrangian tracers with position $z$ solving

$$\frac{dz_l}{dt} = X(z_l, t), \qquad z_l(0) = z_{l,0}, \qquad l = 1, \ldots, L.$$

For simplicity assume that we observe all tracers $z$ at the same set of times $\{t_m\}_{m=1}^M$ and that the $z_{l,0}$'s are known to us:

$$W_{l,m} = z_l(t_m) + \xi_{l,m}, \qquad l = 1, \ldots, L, \ m = 1, \ldots, M,$$

so that, similarly to the Eulerian case, we may write

$$W = \tilde{z} + \xi$$

with $W = (W_{1,1}, W_{2,1}, \ldots, W_{L,M})$, $\tilde{z} = (z_1(t_1), z_2(t_1), \ldots, z_L(t_M))$ and noise $\xi \sim \mathcal{N}(0, \Sigma)$ for some covariance matrix $\Sigma$. Figure 2 illustrates the set-up, showing a snapshot of the flow field streamlines for $X(v, t)$ and the tracer locations $z_l(t)$ at some time instance $t$. We refer to this problem as *Lagrangian data assimilation*, since the observations are of Lagrangian particle trajectories.

Note that, for both Eulerian and Lagrangian observations, $\tilde{z}$ is a function of the initial condition $X_0$ of the Navier–Stokes equation.

In applications to weather forecasting, compressible fluid flow models would in fact be more appropriate. We have chosen to use the incompressible Navier–Stokes equation to model the fluid in both examples simply to unify the presentation of the Eulerian and Lagrangian data models. A more realistic model for weather forecasting would employ a nondissipative model for the velocity field, supporting waves, such as the shallow water equations [4, 19]. The techniques described in this section could be generalized to such models.

5.1. *The Fourier expansion.*  We now construct the probability measure of interest, namely the probability of $X_0$ given $Y$, for both the Eulerian and Lagrangian problems. Any mean-zero $X \in L^2(\mathbb{T}^2, \mathbb{C}^2)$ can be expanded as a Fourier series in the form

$$X(v) = \sum_{k \in \mathbb{Z}^2 \setminus \{0\}} x_k \exp(2\pi i k \cdot v).$$

The divergence-free condition is equivalent to setting $x_k \cdot k = 0$ for all $k$, so we can form an orthonormal basis for $\mathcal{H}$ by letting

$$e_k = \frac{k^\perp}{|k|} \exp(2\pi i k \cdot x),$$

where $k^\perp = (k_2, -k_1)$. Then any $X \in \mathcal{H}$ may be written as

$$X = \sum_{k \in \mathbb{Z}^2 \setminus \{0\}} x_k e_k = \sum_k \langle X, e_k \rangle e_k.$$

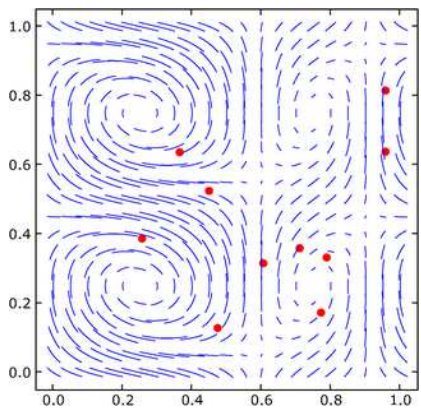

Fig. 2.   *An example configuration of the velocity field at a given time instance. The small circles correspond to a number of Langangian tracers.*

We want to restrict attention to real-valued functions, whence the above expansion can be restated as

$$X = \sum_{k \in \mathbb{Z}_U^2} \{x_k^{\sin} e_k^{\sin} + x_k^{\cos} e_k^{\cos}\} \tag{5.1}$$

for the orthonormal $e_k^{\sin}(u) = \frac{k^\perp}{|k|}\sqrt{2}\sin(2\pi k \cdot u)$, $e_k^{\cos}(u) = \frac{k^\perp}{|k|}\sqrt{2}\cos(2\pi k \cdot u)$ and the upper half of $\mathbb{Z}^2 \setminus \{0\}$,

$$\mathbb{Z}_U^2 = \{k = (k_1, k_2) \in \mathbb{Z}^2 : k_1 > 0, k_2 \geq 0 \text{ or } k_1 \leq 0, k_2 > 0\}.$$

The Stokes operator $A$ is such that

$$Ae_k = 4\pi^2 |k|^2 e_k. \tag{5.2}$$

We will assign a prior measure, $\tilde{\pi}$, on $X_0$. We choose this to be a Gaussian measure with mean zero and precision operator $\mathcal{L} = -A^\alpha$ for a positive real $\alpha$. To be precise, the covariance operator will be

$$\mathcal{C} = \sum_{k \in \mathbb{Z}_U^2} (4\pi^2)^{-\alpha} |k|^{-2\alpha} (e_k^{\sin} \otimes e_k^{\sin} + e_k^{\cos} \otimes e_k^{\cos}).$$

For the Gaussian measure to be well defined it is necessary that $\mathcal{C}$ is trace-class, that is, $\sum_{k \in \mathbb{Z}_U^2} |k|^{-2\alpha} < \infty$, which requires $\alpha > 1$. We condition the prior on the observations, to find the posterior measure, $\pi$, on $X_0$. As noted before, $\tilde{z}$ is a function (the so-called *observation operator*) $\mathcal{G}$ of $X_0$, the initial condition, so we may write

$$W = \mathcal{G}(X_0) + \xi,$$

where one should have in mind that $\mathcal{G} = \mathcal{G}^{\mathrm{EUL}}$ and $\mathcal{G} = \mathcal{G}^{\mathrm{LAG}}$ for the Eulerian and Lagrangian case, respectively. The likelihood of $W$ conditionally on $X_0$ is

$$\mathbb{P}[Y \mid X_0] \propto \exp(-\tfrac{1}{2}|W - \mathcal{G}(X_0)|_\Sigma^2).$$

By the Bayes rule we deduce that

$$\frac{d\pi}{d\tilde{\pi}}(X_0) \propto \exp\left(-\frac{1}{2}|W - \mathcal{G}(X_0)|_\Sigma^2\right), \qquad \tilde{\pi} = \mathcal{N}(0, \mathcal{C}). \tag{5.3}$$

We have thus constructed another explicit example of the structure (3.1) where now

$$\phi(X) = \tfrac{1}{2}|W - \mathcal{G}(X)|_\Sigma^2.$$

By tuning $\alpha$ we may induce more smoothness in the prior and posterior measures. From the Karhunen–Loève expansion, under the prior Gaussian measure $\tilde{\pi}$ we get

$$x_k^{\sin}, x_k^{\cos} \sim \mathcal{N}(0, (4\pi^2)^{-\alpha}|k|^{-2\alpha}), \tag{5.4}$$

all coefficients being independent of each other. A finite-dimensional approximation of (5.3) will be derived by truncating the Fourier series. For integer $N > 0$ we define

$$\mathbb{Z}_{U,N}^2 := \{k \in \mathbb{Z}_U^2 : |k_1|, |k_2| \le N\}.$$

One can check that the cardinality of $\mathbb{Z}_{U,N}^2$ is $2N(N+1)$. We will arrange $x_k^{\sin}, x_k^{\cos}$ with $k \in \mathbb{Z}_{U,N}^2$ into an $n$-dimensional vector $x$, with $n = 4N(N+1)$. To this end, we order the elements of $\mathbb{Z}_U^2$ spiral-wise, as shown in Figure 3; the construction gives rise to an ordering mapping $\sigma : \mathbb{Z}_U^2 \mapsto \mathbb{Z}^+$, analytically specified as

$$\sigma((i,j)) = 2(N-1)N + M,$$

where

$$N = N((i,j)) = |i| \vee j,$$

$$M = M((i,j)) = \begin{cases} j+1, & i = N, j < N, \\ 2N-i+1, & j = N, i > -N, \\ 4N-j+1, & i = -N. \end{cases}$$

The mapping is better understood via Figure 3. We set

$$(5.5) \qquad x = (x_{\sigma^{-1}(i)}^{\sin}, x_{\sigma^{-1}(i)}^{\cos})_{i=1}^{2N(N+1)}.$$

We can now approximate (5.3) as follows:

$$\frac{d\pi_n}{d\bar{\pi}_n}(x) \propto \exp(-\phi(P_N(x))), \qquad P_N(x) := \sum_{k \in \mathbb{Z}_{U,N}^2} \{x_k^{\sin} e_k^{\sin} + x_k^{\cos} e_k^{\cos}\}.$$

(5.6)

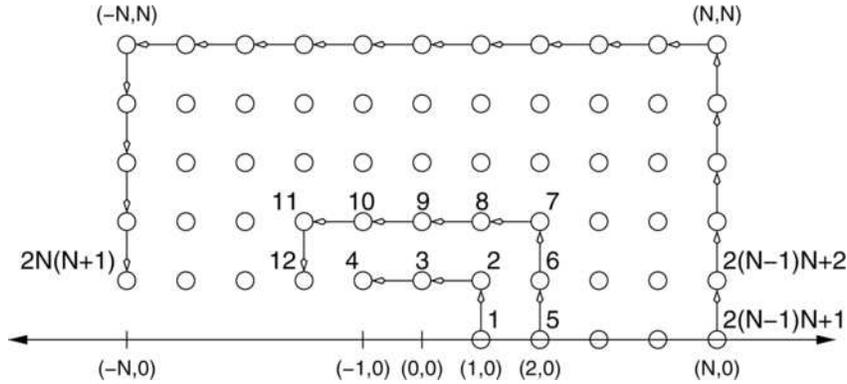

Fig. 3. *The construction of the ordering mapping $\sigma : \mathbb{Z}_U^2 \mapsto \mathbb{Z}^+$. The circles represent points in $\mathbb{Z}_U^2$. The paired numbers in parentheses are coordinates of points; the single numbers show the ordering of points w.r.t. $\sigma(\cdot)$. For example, $\sigma((2,0)) = 5$, $\sigma((2,2)) = 7$.*

5.2. *The results.* In this context, for $s \in \mathbb{R}$ we define

$$\|X\|_s = \left( \sum_{k \in \mathbb{Z}_U^2} |k|^{2s}((x_k^{\sin})^2 + (x_k^{\cos})^2) \right)^{1/2}$$

with corresponding space $\mathcal{H}^s = \{X : \|X\|_s < \infty\}$. We will need some analytical results from [11]. Assume that the force $h$ in the Navier–Stokes equation is such that $\int_0^T \|h\|_\varepsilon^2 \, dt < \infty$ for some $\varepsilon > 0$. Then, in the Eulerian case it is true that

$$|\mathcal{G}^{\mathrm{EUL}}(X)| \leq C(1 + \|X\|^2)$$

for any $X \in \mathcal{H}$. Also, there exists locally bounded function $\delta : \mathbb{R}^+ \times \mathbb{R}^+ \mapsto \mathbb{R}^+$ such that

$$|\mathcal{G}^{\mathrm{EUL}}(Y) - \mathcal{G}^{\mathrm{EUL}}(X)| \leq \delta(\|X\|, \|Y\|)\|Y - X\|.$$

For the Langrangian case we have

$$|\mathcal{G}^{\mathrm{LAG}}(X)| \leq C(1 + \|X\|^2)^{3/2}$$

and, if $X, Y \in \mathcal{H}^s$ for some $s > 0$, then there exists locally bounded $\delta$ such that

$$|\mathcal{G}^{\mathrm{LAG}}(Y) - \mathcal{G}^{\mathrm{LAG}}(X)| \leq \delta(\|X\|_s, \|Y\|_s)\|Y - X\|_s.$$

Note that $E_{\tilde{\pi}}\|X\|_s^2 < \infty$ for any $s < \alpha - 1$, so we can restrict attention on the set of full probability measure $\mathcal{H}^{s_0}$ for some chosen $0 < s_0 < \alpha - 1$ instead of the whole $\mathcal{H}$. Under this observation, one can see that the conditions on $\mathcal{G}^{\mathrm{LAG}}$ and $\mathcal{G}^{\mathrm{EUL}}$ both imply that for $Y, X \in \mathcal{H}^{s_0}$:

$$
\begin{aligned}
(5.7) \qquad |\phi(Y) - \phi(X)| &\leq \delta(\|X\|_{s_0}, \|Y\|_{s_0})\|Y - X\|_{s_0}, \\
|\phi(X)| &\leq C(1 + \|X\|^6)
\end{aligned}
$$

for some locally bounded function $\delta$; these are reminiscent of Condition 5 which, as we show in Appendices A and B, they imply after making the standard correspondence between $X$ and its Fourier coefficients and employing a spectral truncation to obtain a sampling problem in $\mathbb{R}^n$.

PROPOSITION 3. *For any $\alpha > 1$, the limiting results of Theorems 3 and 4 apply to the target $\pi_n$ in (5.6) both for $\mathcal{G} = \mathcal{G}^{\mathrm{EUL}}$ and $\mathcal{G} = \mathcal{G}^{\mathrm{LAG}}$, under the specifications*

$$\kappa = \alpha/2, \qquad \tau^{\mathrm{RWM}} = \frac{1}{2}\pi^{2\alpha} \int_{[0,1]^2} (x^2 + y^2)^\alpha \, dx \, dy,$$

$$\tau^{\mathrm{SLA}} = \frac{1}{2}\pi^{6\alpha} \int_{[0,1]^2} (x^2 + y^2)^{3\alpha} \, dx \, dy.$$

**6. Further directions.** We present some additional MCMC samplers which can give improved results.

6.1. *Preconditioning.* The local MCMC samplers we have so far considered use the same step-size over all directions of the target domain. The standard deviations decrease as $i^{-k}$ and the sampler will adjust its step-size to accommodate for the smallest scale, $n^{-k}$, resulting in the $\mathcal{O}(n^{-2\kappa-1})$ squared-jump-distance reported in the paper. When the different scales are analytically known one can think of allowing the step-size to vary over the different coordinates. We will now briefly pursue such a direction: we demonstrate that the analytical conditions stated so far in the paper suffice to also describe these slightly modified algorithms.

We follow the general context with target distribution $\pi_n$ in (1.3) and scalings $\lambda_i = \lambda_{i,n}$ with $C_- i^{-\kappa} \leq \lambda_{i,n} \leq C_+ i^{-\kappa}$ for some $\kappa \geq 0$. Consider the following MCMC proposals:

$$(6.1) \qquad \mathrm{P-RWM}:\ y_i = x_i + \sigma_n \lambda_{i,n} Z_i,$$

$$(6.2) \qquad \mathrm{P-SLA}:\ y_i = x_i + \frac{\sigma_n^2 \lambda_{i,n}}{2}(-g'(x_i/\lambda_{i,n})) + \sigma_n \lambda_{i,n} Z_i,$$

which compared with the original RWM and SLA, adjust the step-size to the individual scales and use step-size $\sigma_n \lambda_{i,n}$ for the $i$th coordinate instead of $\sigma_n$. The prefix (P-) stands for "preconditioning" referring to the effect of these algorithms to flatten the standard deviations along the scalar components before applying a regular RWM or SLA proposal. For the following result we set $\Lambda = \Lambda_n = \mathrm{diag}\{\lambda_{1,n}, \ldots, \lambda_{n,n}\}$.

COROLLARY 2. *Let $\pi_n$ in (1.3) be the target distribution for P-RWM and P-SLA. Assume that $\sigma_n^2 = l^2 n^{-\rho}$ for $l, \rho > 0$. If $\{\phi_n\}$ satisfies Conditions 3 and 5, then:*

(i) *if $\rho = I$, then $E[a_n(x,y)] \to a(l)$ and*

$$S_n^P := E\left[\left(\frac{x'_{i*} - x_{i*}}{\lambda_{i*,n}}\right)^2\right] = l^2 a(l) \times n^{-I} + o(n^{-I})$$

*for $a(l)$ as in Theorem 1 under the specification $\tau^{\mathrm{RWM}} = \tau^{\mathrm{SLA}} = 1$,*

(ii) *if $\rho > I$, then $E[a_n(x,y)] \to 1$ and $S_n^P = l^2 n^{-\rho} + o(n^{-\rho})$,*

(iii) *if $0 < \rho < I$, then $E[a_n(x,y)] = \mathcal{O}(n^{-p})$ and $S_n^P = \mathcal{O}(n^{-p})$ for any $p \geq 0$.*

PROOF. On the transformed space $x \mapsto \Lambda^{-1}x$ algorithms P-RWM and P-SLA coincide with RWM and SLA, respectively, as applied to the target $\pi_n^t$ specified as

$$\frac{d\pi_n^t}{d\tilde{\pi}_n^t}(x) = \exp\{-\phi_n(\Lambda x)\}, \qquad \tilde{\pi}_n^t(x) = \prod_{i=1}^n f(x_i).$$

It is a straightforward consequence of the definition of the $|x|_s$-norm that if $\phi_n$ satisfies Condition 5, then so does $\phi_n^t(x) = \phi_n(\Lambda x)$. So, all results of Theorems 3 and 4 apply for the target $\pi_n^t$. The statements of the corollary correspond only to a translation of these results on the original coordinates $x = \Lambda x^t$. $\square$

6.2. *Implicit scheme and Gaussian law.* An implicit version of SLA applied when the reference measure is Gaussian satisfies an identity allowing, in some cases, for $\mathcal{O}(1)$ squared-jump-distance. To see this, assume that $\tilde{\pi}_n$ is Gaussian with $f \equiv N(0, 1)$ and consider the proposal $y$ for current position $x$ defined through the equation

$$(6.3) \qquad \theta\text{-SLA: } y_i = x_i + \frac{\sigma_n^2}{2}\left(-\theta\frac{y_i}{\lambda_i^2} - (1-\theta)\frac{x_i}{\lambda_i^2}\right) + \sigma_n Z_i$$

for $\theta \in [0, 1]$. Note that $\theta = 0$ corresponds to the standard SLA proposal (1.5). After some calculations we find that for this proposal the acceptance probability is

$$a_n(x, y) = 1 \wedge \exp\left\{\phi_n(x) - \phi_n(y) + \left(\theta - \frac{1}{2}\right)\frac{\sigma_n^2}{4\lambda_i^2}(x_i^2/\lambda_i^2 - y_i^2/\lambda_i^2)\right\}.$$

Clearly, $a_n(x, y) \equiv 1$ if $\theta = 1/2$ and $\phi_n \equiv 0$, so Gaussian targets are invariant under the update (6.3) with $\theta = 1/2$ (which corresponds to the familiar *trapezoidal rule* in the numerical analysis literature) for any step-size $\sigma_n$. Thus, in the Gaussian scenario, we get $\mathcal{O}(1)$ squared-jump-distance. Even for the case when the change of measure makes the target non-Gaussian, we have shown in [5] and [6] that such implicit proposals can provide well defined MCMC algorithms for infinite-dimensional targets having a density w.r.t. a Gaussian measure thus giving rise to algorithms of $\mathcal{O}(1)$ squared-jump-distance in the discretized context.

**7. Conclusions.** We have presented a family of nonproduct distributions, arising naturally in high-dimensional applications, for which analytical results for the asymptotic behavior of local MCMC algorithms can be obtained. The results in the paper constitute a contribution toward the understanding of the computational complexity of Metropolis–Hastings methods in application to high-dimensional, complex target measures with mathematical structure tailored to a range of important applications.

In this context, the inverse of the squared-jump-distance $S_n$ provides a measure of the number of steps required to cover the invariant measure, as a function of the dimension of the state space, $n$. In a concrete application this needs to be combined with information about the cost of each proposed step, again as a function of $n$. To be concrete we consider the case where the reference measure is a mean-zero Gaussian one with covariance operator $\mathcal{C}_n$ and precision operator $\mathcal{L}_n = -\mathcal{C}_n^{-1}$ and observe that, although our analysis in the paper is conducted in the Karhunen–Loève basis, we need not assume that this is known to implement the methods. That is, we may write the proposals as:

$$\text{RWM: } y = x + \sigma_n Z,$$

$$\text{SLA: } y = x + \frac{\sigma_n^2}{2}\mathcal{L}_n x + \sigma_n Z,$$

$$\text{P-RWM: } y = x + \sigma_n \mathcal{C}_n^{1/2} Z,$$

$$\text{P-SLA: } y = x - \frac{\sigma_n^2}{2}x + \sigma_n \mathcal{C}_n^{1/2} Z,$$

$$\theta\text{-SLA: } y = x + \frac{\sigma_n^2}{2}(\theta\mathcal{L}_n y + (1-\theta)\mathcal{L}_n x) + \sigma_n Z.$$

Having in mind also the calculation of the acceptance probability, all methods require evaluation of $\phi_n(x)$ and $\mathcal{L}_n x$; then, P-RWM and P-SLA require drawing from $N(0, \mathcal{C}_n)$ and $\theta$-SLA (for $\theta = 1/2$) inverting $I - \frac{\sigma_n^2}{4}\mathcal{L}_n$. All such costs should be taken into account for an overall comparison of the different algorithms. Thus, even in the case of Gaussian reference measure, the relative efficiency of the methods depends crucially on the precise structure of the reference measure; for instance, the case of Markovian reference measures, for which the precision operator has a banded structure, will be markedly different from arbitrary non-Markovian problems.

From a mathematical point of view the results in this paper center on a delicate interplay between properties of the reference measure and properties of the change of measure. The tail properties of the reference measure are captured in the scaling properties of the standard deviations [see (1.3c) and (2.6)]. The assumptions we make about $\phi$ in Conditions 3, 4 and 5 control the manner in which the target measure differs from the reference measure, in the tails. Since the tails control the optimal scaling of the algorithms this is key to the analysis. In particular the conditions on the exponents used in the norms in Conditions 4 and 5, and their upper bounds in terms of $\kappa$, are precisely those required to ensure that the product measure dominates in the tails; the choice of norms really matters as we approach an infinite-dimensional limit. It is notable that the conditions imposed on the change of measure do indeed hold for the complex infinite-dimensional applications that we consider in this paper.

We anticipate that there is a whole range of other applications which fall into the class of problems we consider in this paper. For example, one other natural application area is the case when $\tilde{\pi}_n$ represents a prior of independent components for a Bayesian analysis with $d\pi_n/d\tilde{\pi}_n$ corresponding to the likelihood function. In this context, for instance, the lower bound for $\phi_n$ in Condition 3 (resp. upper bound for the likelihood) will be typically satisfied.

A particularly interesting avenue for further study in this area concerns the possibility of obtaining diffusion limits for the MCMC methods when the proposal steps are scaled optimally. This program has been carried out in the case of i.i.d. product targets in [21] and [22] (see also [1] and [2]). In that case each component is asymptotically independent from the others, so it is possible to consider an independent scalar diffusion process in each coordinate. In the nonproduct scenario considered in this paper it is anticipated that the diffusion limit will be a stochastic PDE which is reversible w.r.t. the target measure. Proving such a result would yield further insight into the optimality of MCMC methods.

## APPENDIX A: TAYLOR EXPANSIONS FOR MCMC ALGORITHMS (TARGET $\tilde{\pi}_N$)

We present a Taylor expansion for the term $R_n$ in (2.1) to be used throughout the proofs. In this context we assume that the target is $\tilde{\pi}_n$ with $x \sim \tilde{\pi}_n$.

**A.1. The RWM algorithm.** We consider the exponent $R_n \equiv R_n^{\text{RWM}}$ in (2.2).

CASE A ($\sigma_n^2 = l^2 n^{-2\kappa-1} n^{-\varepsilon}$ with $\varepsilon \geq 0$). Viewing $R_n$ as a function of $\sigma_n$, we employ a second-order Taylor expansion around $\sigma_n = 0$, separately for each index $i = 1, \ldots, n$. Thus,

$$(A.1) \qquad R_n = \mathcal{A}_{1,n} + \mathcal{A}_{2,n} + \mathcal{U}_n,$$

where we have defined

$$\mathcal{A}_{1,n} = \sum_{i=1}^n \frac{\sigma_n}{\lambda_i} C_{1,i}, \qquad \mathcal{A}_{2,n} = \sum_{i=1}^n \frac{\sigma_n^2}{\lambda_i^2} C_{2,i}, \qquad \mathcal{U}_n = \sum_{i=1}^n \frac{\sigma_n^3}{\lambda_i^3} U_{i,n},$$

$$(A.2) \qquad C_{1,i} = -g'(x_i/\lambda_i) Z_i, \qquad C_{2,i} = -g''(x_i/\lambda_i) Z_i^2/2,$$

$$U_{i,n} = -g^{(3)}(x_i/\lambda_i + Z_i \sigma_i^*/\lambda_i) Z_i^3/6$$

for some $\sigma_i^* \in [0, \sigma_n]$ different for each $i$. Note that $\{C_{1,i}\}_i$, $\{C_{2,i}\}_i$ are both sequences of i.i.d. random variables. Using Condition 1(i) we find that

$$(A.3) \qquad |U_{i,n}| \leq M_1(x_i/\lambda_i) M_2(Z_i) M_3(\sigma_i^*/\lambda_i)$$

for some positive polynomials $M_1$, $M_2$, $M_3$. From Condition 1(ii) we get $E[M_1(x_i/\lambda_i)] < \infty$; also $E[M_2(Z_i)] < \infty$, both these expectations not depending on $i$. Note that $\sigma_n/\lambda_i \to 0$, so $M_3(\sigma_i^*/\lambda_i) \leq K_0$ for a constant $K_0 > 0$. We can now obtain the following results:

- $\mathcal{U}_n \xrightarrow{L_2(\tilde{\pi}_n)} 0$,
- if $\varepsilon > 0$, then $\mathcal{A}_{1,n} \xrightarrow{L_2(\tilde{\pi}_n)} 0$, $\mathcal{A}_{2,n} \xrightarrow{L_2(\tilde{\pi}_n)} 0$,
- if $\varepsilon = 0$, then $\mathcal{A}_{1,n} \xrightarrow{\mathcal{L}} N(0, l^2\frac{K}{1+2\kappa})$, $\mathcal{A}_{2,n} \xrightarrow{L_2(\tilde{\pi}_n)} -\frac{l^2}{2}\frac{K}{1+2\kappa}$, where

(A.4) $$K = K^{\mathrm{RWM}} = E_f[(g'(X))^2].$$

The limit for $\mathcal{U}_n$ follows directly from (A.3). Also, the stated limits for $\varepsilon > 0$ follow from simple calculations. For the results when $\varepsilon = 0$ we note that a version of the Lindeberg–Feller theorem for the case of scaled sums of i.i.d. variables (see, e.g., Theorem 2.73 of [7]) gives that $\mathcal{A}_{1,n}$ converges in law to a zero-mean Gaussian random variable with variance the limit as $n \to \infty$ of $E[\mathcal{A}_{1,n}^2]$ which can be easily found to be $l^2 E[C_{1,\cdot}^2]/(1+2\kappa)$. For $\mathcal{A}_{2,n}$, straightforward calculations give that it converges in $L_2(\tilde{\pi}_n)$ to $l^2 E[C_{2,\cdot}]/(1+2\kappa)$; the product law gives that

$$E_f[g'(X)^2] = E_f[g''(X)],$$

so, in fact, the limit for $\mathcal{A}_{2,n}$ can also be expressed in terms of $K$.

CASE B ($\sigma_n^2 = l^2 n^{-2\kappa-1} n^\varepsilon$ with $\varepsilon \in (0,1)$). We now select a positive integer $m$ such that

$$m + 1 > 2/(1-\varepsilon)$$

and use the $m$-order expansion:

$$R_n = \sum_{j=1}^m \mathcal{A}_{j,n} + \mathcal{U}'_n$$

in the place of (A.1), where now

$$\mathcal{A}_{j,n} = \sum_{i=1}^n \frac{\sigma_n^j}{\lambda_i^j} C_{j,i}, \qquad \mathcal{U}'_n = \sum_{i=1}^n \frac{\sigma_n^{m+1}}{\lambda_i^{m+1}} U'_{i,n},$$

$$C_{j,i} = -g^{(j)}(x_i/\lambda_i) Z^j/(j!), \qquad U'_{i,n} = -g^{(m+1)}(x_i/\lambda_i + Z_i \sigma_i^*/\lambda_i)\frac{Z_i^{m+1}}{(m+1)!}$$

for some corresponding $\sigma_i^* \in [0, \sigma_n]$. We can now obtain the following results:

- $\mathcal{U}'_n \xrightarrow{L_1(\tilde{\pi}_n)} 0$,
- $E[R_n] \to -\infty$ as fast as $-n^\varepsilon$,
- for any positive integer $q$, $E[(R_n - E[R_n])^{2q}] = \mathcal{O}(n^{\varepsilon \cdot q})$.

For the first result note that the residual terms $U'_{i,n}$ can be bounded as in (A.3), so the limit follows from the particular choice of $m$. For $E[R_n]$ we remark that, since $E[\mathcal{A}_{1,n}] = 0$, we have

$$E[R_n] = \sum_{j=2}^m E[\mathcal{A}_{j,n}] + \mathcal{O}(1), \qquad E[\mathcal{A}_{j,n}] = \sum_{i=1}^n \frac{\sigma_n^j}{\lambda_i^j} E[C_{j,\cdot}].$$

From the analytical expression for $C_{2,i}$:

$$E[C_{2,\cdot}] = -\frac{1}{2} \int_{\mathbb{R}} \{g'(u)\}^2 \exp\{-g(u)\} \, du < 0.$$

All other $C_{j,\cdot}$ satisfy $E|C_{j,\cdot}| < \infty$. Trivially, the highest-order term among the summands for $E[R_n]$ is the one corresponding to $j = 2$ which indeed grows to $-\infty$ as fast as $-n^\varepsilon$. For the third result, among the $(m+1)$ zero-mean sums comprising $R_n - E[R_n]$ the one with the highest-order $L_{2q}$-norm is $\mathcal{A}_{1,n}$, so the triangle inequality gives that the provided expectation is of the same order as $E[\mathcal{A}_{1,n}^{2q}]$. Now we can take out of the expectation the $n^{\varepsilon/2}$ factor of $\sigma_n$; the remaining expectation is $\mathcal{O}(1)$. To prove this last statement one needs to consider the polynomial expansion; we avoid further details.

**A.2. The SLA algorithm.** We now use Taylor expansions for the corresponding term $R_n \equiv R_n^{\mathrm{SLA}}$ given in (2.3).

CASE A ($\sigma_n^2 = l^2 n^{-2\kappa - 1/3} n^{-\varepsilon}$ with $\varepsilon \geq 0$). We employ a sixth-order Taylor expansion to obtain the structure

(A.5)
$$R_n = \sum_{j=3}^{6} \mathcal{A}_{j,n} + \mathcal{U}_n,$$

$$\mathcal{A}_{j,n} = \sum_{i=1}^{n} \frac{\sigma_n^j}{\lambda_i^j} C_{j,i}, \qquad \mathcal{U}_n = \sum_{i=1}^{n} \frac{\sigma_n^7}{\lambda_i^7} U_{i,n}.$$

We start the summation from $j = 3$ since the first two sums are identically zero. The above expansion considers some corresponding i.i.d. sequences $\{C_{j,i}\}_{i=1}^{n}$ for each $j$ and residual terms; any $C_{j,i}$ is a polynomial function (the same over $i$) of $Z_i$ and $x_i/\lambda_i$. Using the calculations in [22] we have

$$C_{1,\cdot} = C_{2,\cdot} \equiv 0, \qquad E[C_{3,\cdot}] = E[C_{4,\cdot}] = E[C_{5,\cdot}] = 0, \qquad E[C_{6,\cdot}] < 0.$$

From this paper we will also require the analytical calculation

$$C_{3,i} = g^{(3)}(x_i/\lambda_i) Z_i^3/12 - g'(x_i/\lambda_i) g''(x_i/\lambda_i) Z_i/4.$$

From the product rule,

$$E_f[(g'(X)g''(X))^2] - 2E_f[(g'g''g^{(3)})(X)] = E_f[g''(X)^3],$$

thus

$$E[C_{3,\cdot}^2] = E_f[3g''(X)^3 + 5(g^{(3)}(X))^2]/48.$$

We will need the following results:

- $\mathcal{U}_n \xrightarrow{L_2(\tilde{\pi}_n)} 0$,

- if $\varepsilon > 0$, then $\mathcal{A}_{j,n} \xrightarrow{L_2(\tilde{\pi}_n)} 0$ for all $3 \le j \le 6$,
- if $\varepsilon = 0$, then $\mathcal{A}_{3,n} \xrightarrow{\mathcal{L}} N(0, l^6 K)$, both $\mathcal{A}_{4,n}, \mathcal{A}_{5,n} \xrightarrow{L_2(\tilde{\pi}_n)} 0$ and $\mathcal{A}_{6,n} \xrightarrow{L_2(\tilde{\pi}_n)} -\frac{l^6}{2} \frac{K}{1+6\kappa}$, where

(A.6) $$K = K^{\mathrm{SLA}} = E_f[3g''(X)^3 + 5(g^{(3)}(X))^2]/48 > 0.$$

The residuals can be bounded as in (A.3); the limit is then straightforward. For the other quantities we work as in the case of RWM; note that the identity

$$E[C_{3,\cdot}^2] + 2E[C_{6,\cdot}] = 0$$

demonstrated in [22] allows for the limits of $\mathcal{A}_{3,n}$ and $\mathcal{A}_{6,n}$ to be expressed in terms of the same constant $K$.

CASE B ($\sigma_n^2 = l^2 n^{-2\kappa - 1/3} n^\varepsilon$ with $\varepsilon \in (0, 1/3)$). We now use an $m$-order expansion, where $m$ is such that

$$(m+1)(1 - 3\varepsilon) \ge 6$$

and consider the corresponding expansion

$$R_n = \sum_{j=1}^m \mathcal{A}_{j,n} + \mathcal{U}_n'.$$

Working as for Case B of RWM we obtain the following results:

- $\mathcal{U}_n' \xrightarrow{L_1(\tilde{\pi}_n)} 0$,
- $E[R_n] \to -\infty$ as fast as $n^{3\varepsilon}$,
- for any positive integer $q$, $E[(R_n - E[R_n])^{2q}] = \mathcal{O}(n^{3\varepsilon \cdot q})$.

## APPENDIX B: PROOFS

The following lemmas will be used at the proofs.

LEMMA B.1.   *Let $T$ be a random variable. Then:*

(i) *for any $\gamma > 0$:*

$$E[1 \wedge e^T] \ge e^{-\gamma} \left(1 - \frac{E|T|}{\gamma}\right),$$

(ii) *if $E[T] = -c < 0$, then for the residual $T_{\mathrm{res}} = T - E[T]$ and any $q > 0$:*

$$E[1 \wedge e^T] \le e^{-c/2} + 2^q \frac{E[|T_{\mathrm{res}}|^q]}{c^q}.$$

PROOF. For the lower bound we have that

$$E[1 \wedge e^T] \geq E[(1 \wedge e^T) \cdot I\{|T| \leq \gamma\}] \geq e^{-\gamma} P[|T| \leq \gamma],$$

which, from the Markov inequality, gives the required result. For the second inequality:

$$E[1 \wedge e^T] = E\left[(1 \wedge e^T) \cdot I\left\{|T_{\mathrm{res}}| \leq \frac{c}{2}\right\}\right] + E\left[(1 \wedge e^T) \cdot I\left\{|T_{\mathrm{res}}| > \frac{c}{2}\right\}\right]$$

$$\leq e^{-c/2} + P\left[|T_{\mathrm{res}}| > \frac{c}{2}\right],$$

which, again from Markov inequality, gives the required result. $\square$

LEMMA B.2. If $X \sim \mathcal{N}(\mu, \sigma^2)$, then

$$E[1 \wedge e^X] = \Phi(\mu/\sigma) + e^{\mu+\sigma^2/2}\Phi(-\sigma - \mu/\sigma).$$

LEMMA B.3 (Sobolev embedding). (i) Let $X \in L^2([0,T], \mathbb{R})$. Consider the expansion w.r.t. the sinusoidal basis $\{\sin(k\pi \cdot /T)\}_{k=1}^{\infty}$:

$$X(t) = \sum_{k=1}^{\infty} x_k \sin(k\pi t/T).$$

If $s, p \in \mathbb{R}$ are such that $s < \frac{1}{2}$ and $2 \leq p < (\frac{1}{2} - s)^{-1}$, then

$$\|X\|_{L_p} \leq C|x|_s$$

for a constant $C > 0$, where $|x|_s = (\sum_{k=1}^{\infty} k^{2s}|x_k|^2)^{1/2}$.

(ii) Let $X = (X_t)_{t=1}^N \in \mathbb{R}^N$ for integer $N > 0$. Consider the sinusoidal basis in $\mathbb{R}^N$, $\{\sin(\frac{k\pi t}{N+1}); t = 1, \ldots, N\}_{k=1}^N$, and the expansion:

$$X_t = \sum_{k=1}^{n} x_k \sin\left(\frac{k\pi t}{N+1}\right).$$

If $s, p \in \mathbb{R}$ are such that $s < \frac{1}{2}$ and $2 \leq p < (\frac{1}{2} - s)^{-1}$, then

$$\left(\sum_{t=1}^{N} |X_t|^p \frac{1}{N+1}\right)^{1/p} \leq C|x|_s$$

for a constant (independent of $X$ and $n$) $C > 0$.

PROOF. (i) If $Y \in L^2([-T, T], \mathbb{C})$ is periodic with period $2T$, its Fourier expansion is

(B.1) $$Y(t) = \sum_{k=-\infty}^{\infty} y_k e^{ik\pi t/T}, \qquad y_k = \frac{1}{2T}\langle Y, e^{ik\pi \cdot /T}\rangle.$$

The Sobolev embedding (see page 436 of [23]) gives that, for $s, p$ as in the statement of the lemma,

$$(\text{B.2}) \qquad \|Y\|_{L_p} \leq C \left\{ \sum_{k=-\infty}^{k=\infty} (1 + |k|^{2s}) |y_k|^2 \right\}^{1/2}.$$

This is a consequence of the fact that for conjugate positive reals $p$, $q$ (i.e., $p^{-1} + q^{-1} = 1$) we have

$$(\text{B.3}) \qquad \|Y\|_{L_p} \leq C \|y\|_{l_q}.$$

See the above reference for more details. Assume now that $Y$ is specified as follows: $Y(t) = X(t)$ for $t \geq 0$, and $Y(t) = -X(-t)$ when $t < 0$. This symmetricity around the origin means that $y_{-k} = -y_k$ (for instance, $y_0 = 0$), so using this:

$$Y(t) = \sum_{k=1}^{\infty} 2iy_k \sin(k\pi t/T).$$

$X$, $Y$ coincide on $[0, T]$ so $x_k = 2iy_k$ and, using (B.2):

$$\|X\|_{L_p([0,T],\mathbb{R})} \leq \|Y\|_{L_p([-T,T],\mathbb{R})} \leq C|x|_s.$$

(ii) Consider a vector $Y = \{Y_t\}_{t=-(N+1)}^{N} \in \mathbb{C}^{2N+2}$. We define a function $\tilde{Y}$ on $[-T, T]$:

$$\tilde{Y}(t) = \sum_{j=-(N+1)}^{N} Y_j \mathbb{I}_{[jT/(N+1),(j+1)T/(N+1))}(t), \qquad \tilde{Y}(T) = Y_{-(N+1)}.$$

The continuous- and discrete-time Fourier expansions can be written as follows:

$$\tilde{Y}(t) = \sum_{k=-\infty}^{\infty} \tilde{y}_k e^{ik\pi t/T}, \qquad Y_t = \sum_{k=-(N+1)}^{N} y_k e^{ik\pi t/(N+1)},$$

where, after some calculations, we find that for $k \in \mathbb{Z}$ with $k \neq 0$:

$$\tilde{y}_k = -\frac{1}{2ik\pi} \sum_{j=-(N+1)}^{N} Y_j (e^{-ik\pi(j+1)/(N+1)} - e^{-ik\pi j/(N+1)})$$

with $\tilde{y}_0 = \frac{1}{2(N+1)} \sum_{j=-(N+1)}^{N} Y_j$. Also, for $-(N+1) \leq k \leq N$:

$$y_k = \frac{1}{2(N+1)} \sum_{j=-(N+1)}^{N} Y_j e^{-ik\pi j/(N+1)}.$$

One can easily now check that, for $-(N+1) \leq k \leq N$,

$$\tilde{y}_k = -\frac{1}{i\pi} \frac{e^{-i\pi k/(N+1)} - 1}{k/(N+1)} y_k,$$

whereas for $M \in \mathbb{Z}$ and $-(N+1) \leq k \leq N$,

$$\tilde{y}_{M(2N+2)+k} = \frac{k}{M(2N+2)+k} \tilde{y}_k.$$

Note that $|e^{i\pi a} - 1| \leq C|a|$ for $a \in [0,1]$. So, using the last two equations and (B.3), for conjugate $p, q$ we get

$$\left( \sum_{t=-(N+1)}^{N} |Y_t|^p \frac{1}{N+1} \right)^{1/p} \equiv \|\tilde{Y}\|_{L_p([-T,T],\mathbb{R})} \leq C \left( \sum_{k=-(N+1)}^{N} |y_k|^q \right)^{1/q} \equiv \|y\|_{l_q}.$$

An application of Holder's inequality gives (see page 437 of [23] for details) that

$$\|y\|_{l_q} \leq C \left\{ \sum_{k=-(N+1)}^{k=N} (1 + |k|^{2s}) |y_k|^2 \right\}^{1/2}$$

for a constant $C$ independent of $\{y_k\}$ and $n$. So, in total:

$$\left( \sum_{t=-(N+1)}^{N} |Y_t|^p \frac{1}{N+1} \right)^{1/p} \leq C \left\{ \sum_{k=-(N+1)}^{k=N} (1 + |k|^{2s}) |y_k|^2 \right\}^{1/2}.$$

To prove the statement (ii) we use the standard method as in (i): we specify the vector $Y = \{Y_t\}$ as $Y_t = X_t$, for $t = 1, \ldots, N$, $Y_0 = 0$, and $Y_t = -X_{-t}$, for $t = -N, \ldots, -1$, $Y_{-(N+1)} = 0$. Then one can find that $y_0 = y_{-(N+1)} = 0$ and $x_k = 2iy_k$ for $k = 1, \ldots, n$. The required result then follows directly.  □

PROOF OF THEOREM 1.  Expectations are in stationarity, so $x \sim \tilde{\pi}_n$ and $y$ is determined from the transitions (1.4) or (1.5) for RWM or SLA, respectively.

*The RWM algorithm.*  Case (i): $\sigma_n^2 = l^2 n^{-2\kappa-1}$.
From Appendix A.1 we get that

$$R_n \xrightarrow{\mathcal{L}} N\left( -\frac{l^2}{2} \frac{K}{1+2\kappa}, l^2 \frac{K}{1+2\kappa} \right).$$

So, Lemma B.2 gives that

$$\lim_{n \to \infty} E[a_n(x,y)] = 2\Phi\left( -\frac{l}{2} \sqrt{\frac{K}{1+2\kappa}} \right).$$

Case (ii): $\sigma_n^2 = l^2 n^{-2\kappa-1} n^{-\varepsilon}$ for $\varepsilon > 0$.

From Appendix A.1 we get that $R_n \to 0$ in $L_1(\tilde{\pi}_n)$. The result follows from the Lipschitz continuity of $x \mapsto 1 \wedge e^x$.

Case (iii): $\sigma_n^2 = l^2 n^{-2\kappa-1} n^{\varepsilon}$ for $\varepsilon \in (0,1)$.

Appendix A.1 gives that $E[R_n] \to -\infty$ as fast as $-n^{\varepsilon}$, and that for arbitrary integer $q > 0$, $E[(R_n - E[R_n])^{2q}] = \mathcal{O}(n^{\varepsilon \cdot q})$. So, Lemma B.1(ii) implies that $E[a_n(x,y)] \to 0$ faster than any polynomial-order.

*The SLA algorithm.* The proof is similar and follows from the results in Appendix A.2. $\square$

PROOF OF THEOREM 2. Here $a_n(x,y) = 1 \wedge e^{\phi_n(x) - \phi_n(y) + R_n}$, for $R_n$ as in the product case. Note now that

$$E_{\pi_n}[1 \wedge e^{\phi_n(x) - \phi_n(y) + R_n}] \leq C E_{\tilde{\pi}_n}[1 \wedge e^{R_n}]$$

for a constant $C > 0$. So, the required result for $2\kappa < \rho < 2\kappa + I$ follows now directly from Theorem 1(iii). For the case when $\rho \geq 2\kappa + I$ note that, using the Taylor expansions for $R_n$ in Appendix A, we can easily find that

$$\limsup_{n \to \infty} E_{\tilde{\pi}_n} |R_n| < \infty.$$

The boundedness condition on $\phi_n$ gives that

$$\limsup_{n \to \infty} E_{\pi_n} |\phi_n(x) - \phi_n(y) + R_n| \leq C_1 + C_2 \limsup_{n \to \infty} E_{\tilde{\pi}_n} |R_n| < \infty$$

for constants $C_1, C_2 > 0$. So, Lemma B.1(i) implies that $E_{\pi_n}[a_n(x,y)]$ is lower bounded. $\square$

PROOF OF THEOREM 3. We will first show that, for any $\rho > 2\kappa$, Conditions 3 and 5 imply that $\phi_n(y) - \phi_n(x) \to 0$ in $L_q(\pi_n)$ for any $q > 0$, for both RWM and SLA. We will then proceed with some calculations to obtain the required results.

Let $\sigma_n^2 = l^2 n^{-2\kappa-\varepsilon}$ for some $\varepsilon > 0$. Recall that under $\tilde{\pi}_n$, $x_i \sim \mathcal{N}(0, i^{-2\kappa})$ independently over $i$. Note that for an arbitrary $s < \kappa - 1/2$,

$$E_{\tilde{\pi}_n} |x|_s^2 = E_{\tilde{\pi}_n} \left[ \sum_{i=1}^n i^{2s} x_i^2 \right] \sim \sum_{i=1}^n i^{2(s-\kappa)} < C$$

for a constant $C > 0$. Similarly, we find that if $q$ is a positive integer, then $E_{\tilde{\pi}_n} |x|_s^{2q} \sim (\sum_{i=1}^n i^{2(s-\kappa)})^q$, so for any integer $q > 0$:

$$(B.4) \qquad E_{\tilde{\pi}_n} |x|_s^{2q} < C.$$

This result directly implies that, for both the RWM and SLA proposal $y$,

$$(B.5) \qquad E_{\tilde{\pi}_n} |y - x|_s^{2q} \to 0.$$

To see this note that, from the triangle inequality, applied for the $|\cdot|_s$-norm, and the definition of the proposal $y$,

$$E_{\tilde{\pi}_n}|y-x|_s^{2q} \leq C(E|\sigma_n Z|_s^{2q} + E_{\tilde{\pi}_n}|z|_s^{2q})$$

for vector $z$ with $z_i \equiv 0$ in the case of RWM, and $z_i = \frac{\sigma_n^2}{2} g'(x_i/\lambda_i)/\lambda_i$ in the case of SLA. Note now that we can write $x_i^2 = i^{-2\kappa}\xi_i^2$ for some i.i.d. random variables $\xi_i^2$ (in this case $\xi_1^2 \sim \chi_1^2$, but the particular distribution is not important for our argument, as long as it has finite moments). Similarly, we can write $(\sigma_n Z_i)^2 = n^{-\varepsilon} n^{-2\kappa}\xi_i^2$ and $z_i^2 = n^{-2\varepsilon}n^{-2\kappa}(n/i)^{-2\kappa}\xi_i^2$ for some i.i.d. positive random variables $\xi_i^2$, different for each case. It is now clear that

$$E|\sigma_n Z|_s^{2q} \leq Cn^{-\varepsilon q}E_{\tilde{\pi}_n}|x|_s^{2q}, \qquad E_{\tilde{\pi}_n}|z|_s^{2q} \leq Cn^{-2\varepsilon q}E_{\tilde{\pi}_n}|x|_s^{2q},$$

which explain (B.5). Given (B.4), (B.5), the triangular inequality implies

$$(B.6) \qquad E_{\tilde{\pi}_n}|y|_s^{2q} < C.$$

We now set $\Delta\phi_n := \phi_n(y) - \phi_n(x)$ and proceed as follows: for any $R > 0$,

$$E|\Delta\phi_n|^q = E[|\Delta\phi_n|^q\mathbb{I}[|x|_s \leq R, |y|_s \leq R]] + E[|\Delta\phi_n|^q\mathbb{I}[|x|_s > R \text{ or } |y|_s > R]]$$
$$\leq \gamma(R)E|y-x|_{s'}^q + CE[(1+|x|_{s''}^{pq}+|y|_{s''}^{pq})\mathbb{I}[|x|_s > R \text{ or } |y|_s > R]]$$
$$\leq \gamma(R)E|y-x|_{s'}^q$$
$$\quad + C(E[1+|x|_{s''}^{2pq}+|y|_{s''}^{2pq}])^{1/2}(\mathbb{P}[|x|_s > R] + \mathbb{P}[|y|_s > R])^{1/2},$$

where $\gamma(R) = \sup_{a \leq R, b \leq R} \delta^q(a,b)$. Let $\varepsilon > 0$. From (B.4) and (B.6) and the Markov inequality, we can choose some $R = R(\varepsilon)$ so that the second term on the right-hand side of the last inequality is smaller than $\varepsilon/2$. Also, (B.5) implies that the first term is smaller than $\varepsilon/2$ for sufficiently large $n$. Thus, for any $q > 0$,

$$\phi_n(y) - \phi_n(x) \xrightarrow{L_q(\tilde{\pi}_n)} 0.$$

Condition 3 gives also that, for any $q > 0$,

$$(B.7) \qquad \phi_n(y) - \phi_n(x) \xrightarrow{L_q(\pi_n)} 0.$$

From the Lipschitz continuity of $x \mapsto 1 \wedge e^x$, for any $\rho > 2\kappa$,

$$(B.8) \qquad E_{\pi_n}[1 \wedge e^{\phi_n(x)-\phi_n(y)+R_n}] - E_{\pi_n}[1 \wedge e^{R_n}] \to 0.$$

We now distinguish between RWM and SLA and the various step-sizes.

*The RWM algorithm.*   We use the expansion $R_n = R_n^{\mathrm{RWM}} = \mathcal{A}_{1,n} + \mathcal{A}_{2,n} + \mathcal{U}_n$ in Appendix A.1.

Case (i): $\sigma_n^2 = l^2 n^{-2\kappa-1}$.

For this step-size we have shown in Appendix A.1 that $\mathcal{A}_{2,n} \to -\frac{l^2}{2}\frac{K}{1+2\kappa}$ and $\mathcal{U}_n \to 0$ in $L_1(\tilde{\pi}_n)$. Due to Condition 3, the same limits will also hold in $L_1(\pi_n)$. Recall from (A.2) that $\mathcal{A}_{1,n} = \sigma_n \sum_{i=1}^n \xi_{i,n}$ for $\xi_{i,n} = -\lambda_i^{-1} g'(x_i/\lambda_i) Z_i$. Due to $Z_i$, for each $n$ the process $\{S_{i,n}\}_{i=1}^n$ with $S_{i,n} = \sigma_n \sum_{j=1}^i \xi_{i,n}$ is a martingale. Under $\tilde{\pi}_n$, the independence among $\xi_{i,n}$ together with some tedious calculations give that

$$\sigma_n^2 \sum_{i=1}^n \xi_{i,n}^2 \overset{L_1(\tilde{\pi}_n)}{\longrightarrow} l^2 \frac{K}{1+2\kappa}.$$

From Condition 3, the same limit holds in $L_1(\pi_n)$. The Martingale CLT from page 58 of [17] now gives that, under $\pi_n$:

$$\mathcal{A}_{1,n} \overset{\mathcal{L}}{\to} N\left(0, l^2 \frac{K}{1+2\kappa}\right).$$

So, comparing with the results for the product case in Appendix A.1, $R_n$ has the same limiting behavior under $\pi_n$ and $\tilde{\pi}_n$, implying that

$$(B.9) \qquad E_{\pi_n}[1 \wedge e^{R_n}] \to a(l)$$

for the same $a(l)$ as in the case when the target law is $\tilde{\pi}_n$. Equations (B.8) and (B.9) give the required result.

Case (ii): $\sigma_n^2 = l^2 n^{-2\kappa-1} n^{-\varepsilon}$ for $\varepsilon > 0$.

For this step-size Appendix A.1 gives that $R_n \to 0$ in $L_1(\tilde{\pi}_n)$ and Condition 3 implies that the same limit holds also in $L_1(\pi_n)$. Equation (B.8) now provides the required result.

Case (iii): $\sigma_n^2 = l^2 n^{-2\kappa-1} n^{\varepsilon}$ for $\varepsilon \in (0,1)$.

From Condition 3, it suffices to show that $n^p E_{\tilde{\pi}_n}[1 \wedge e^{\phi_n(x) - \phi_n(y) + R_n}] \to 0$. In Appendix A.1 we have shown that $E_{\tilde{\pi}_n}[R_n] \to -\infty$ as fast as $-n^{\varepsilon}$; also, for any integer $q > 0$, we have shown that $E_{\tilde{\pi}_n}[(R_n - E_{\tilde{\pi}_n}[R_n])^{2q}] = \mathcal{O}(n^{\varepsilon \cdot q})$. From (B.7), the same orders persist if we replace $R_n$ with $\phi_n(x) - \phi_n(y) + R_n$. The result now follows from Lemma B.1(ii).

*The SLA algorithm.*   We use the corresponding expansion $R_n = \mathcal{A}_{3,n} + \mathcal{A}_{4,n} + \mathcal{A}_{5,n} + \mathcal{A}_{6,n} + \mathcal{U}_n'$. The results for cases (ii) and (iii) are obtained exactly as in the case of RWM. For case (i), the Martingale CLT gives (as for RWM) that $R_n$ has the same limiting behavior under both $\pi_n$ and $\tilde{\pi}_n$, and the required result then follows from (B.8). We avoid further details. □

PROOF OF THEOREM 4.   Note first that

$$(B.10) \qquad S_n = E[(x'_{i^*} - x_{i^*})^2] = E[(y_{i^*} - x_{i^*})^2 1 \wedge e^{\phi_n(x) - \phi_n(y) + R_n}].$$

Simple calculations give that, for any $\rho > 2\kappa$,

$$(B.11) \qquad n^\rho E[(y_{i^*} - x_{i^*})^2] \to l^2, \qquad n^\rho E[(y_{i^*} - x_{i^*})^4]^{1/2} \to l^2.$$

Since $n^\rho E[(y_{i^*} - x_{i^*})^4]^{1/2}$ is $n$-uniformly bounded, the Lipschitz continuity of $x \mapsto 1 \wedge e^x$ and Cauchy–Schwarz inequality imply that any term appearing in the exponential in (B.10) can be replaced with its $L_2(\pi_n)$-limit when considering the limiting behavior of $n^\rho S_n$. To be more precise, we have shown, for instance, in the proof of Theorem 3 that $\phi_n(y) - \phi_n(x) \to 0$ in $L_2(\pi_n)$ for any $\rho > 2\kappa$. This gives

$$|n^\rho S_n - n^\rho E[(y_{i^*} - x_{i^*})^2 1 \wedge e^{R_n}]|$$
$$\leq C n^\rho E[(y_{i^*} - x_{i^*})^4]^{1/2} (E|\phi_n(y) - \phi_n(x)|^2)^{1/2} \to 0.$$

For case (ii) of the proposition, we have shown in Appendix A that $R_n \to 0$ in $L_2(\tilde{\pi}_n)$, so also in $L_2(\pi_n)$ from Condition 3. Now, ignoring $\phi_n(x) - \phi_n(y)$ and $R_n$ from the expression for $S_n$ we get that $n^\rho S_n \to l^2$. For case (iii), we use Cauchy–Schwarz to get

$$n^\rho S_n \leq n^\rho E[(y_{i^*} - x_{i^*})^4]^{1/2} E[a_n(x, y)]^{1/2},$$

so the result follows from (B.11) and Theorem 3(iii). We now focus on the case $\sigma_n^2 = l^2 n^{-2\kappa - I}$ and RWM (whence $I = 1$); the proof for SLA is similar. We use the expansion $R_n = \mathcal{A}_{1,n} + \mathcal{A}_{2,n} + \mathcal{U}_n$ in Appendix A.1. Let $\mathcal{A}_{1,n}^*$, $\mathcal{A}_{2,n}^*$, $\mathcal{U}_n^*$ be the variables derived by omitting the $i^*$th summand from the expansions for $\mathcal{A}_{1,n}$, $\mathcal{A}_{2,n}$, $\mathcal{U}_n$, respectively. We define $R_n^*$ as the sum of these terms. From the analytical expressions in Appendix A.1, it is clear that $R_n - R_n^* \to 0$ in $L_2(\pi_n)$. Thus:

$$n^{2\kappa + 1} S_n - n^{2\kappa + 1} E[(y_{i^*} - x_{i^*})^2 1 \wedge e^{R_n^*}] \to 0.$$

We can now factorize:

$$n^{2\kappa + 1} E[(y_{i^*} - x_{i^*})^2 1 \wedge e^{R_n^*}] = l^2 E[Z_{i^*}^2] \times E[1 \wedge e^{R_n^*}].$$

From the proof of Theorem 3 the last expectation, however, converges to $a(l)$ and the required result is established. □

PROOF OF COROLLARY 1. One needs only to replace $\lambda_i$ with $\lambda_{i,n}$ in all statements of Appendix A and the proofs of Theorems 1–4. When the constants $K^{\text{RWM}}$, $K^{\text{SLA}}$ appear in these proofs (where $\lambda_i = i^{-\kappa}$) they are always divided with $(1 + 2\kappa)$, $(1 + 6\kappa)$, respectively; these values arise as the limits of $n^{-(2\kappa + 1)} \sum_{i=1}^n \lambda_i^{-2}$ and $n^{-(6\kappa + 1)} \sum_{i=1}^n \lambda_i^{-6}$, respectively. Revisiting the proofs shows immediately that in the extended setting of the corollary one should now use $\lim_n n^{-(2\kappa + 1)} \sum_{i=1}^n \lambda_{i,n}^{-2}$ and $\lim_n n^{-(6\kappa + 1)} \sum_{i=1}^n \lambda_{i,n}^{-6}$ in the place of the above limits. □

Proof of Proposition 1. For $x, y \in \mathbb{R}^n$, let $X_N = \sum_{i=1}^N x_{\cdot,i} e_i$ and $Y_N = \sum_{i=1}^N x_{\cdot,i} e_i$. Then, from (4.8) and the equivalence between the norms $\|X_N\|_s$ and $|x|_s$, the latter defined in (2.5), we have

$$
\begin{aligned}
|\phi_n(y) - \phi_n(x)| &= |\phi(Y_N) - \phi(X_N)| \\
&\leq C(1 + \|X_N\|_s^p + \|Y_N - X_N\|_s^p)\|Y_N - X_N\|_{L_2} \\
&\leq C(1 + |x|_s^p + |y - x|_s^p)|y - x|
\end{aligned}
$$

for $s = s(p) < 1/2$. In this context $\kappa = 1$ since the reference measure $\tilde{\pi}_n$ is of the structure $\prod_{i=1}^n \mathcal{N}(0, \Lambda_i^2)$ with

$$
\Lambda_i = \sqrt{\frac{2}{\beta}} \frac{T}{\pi} \left\lceil \frac{i}{d} \right\rceil^{-1},
$$

so clearly $C_- i^{-1} \leq \Lambda_i \leq C_+ i^{-1}$ for appropriate $C_-, C_+ > 0$. Thus, we have found that $\phi_n$ satisfies Condition 4. Then, for the case of RWM:

$$
n^{-2\kappa - 1} \sum_{i=1}^n \Lambda_i^{-2} = n^{-3} \frac{\beta \pi^2}{2T^2} \sum_{i=1}^n \left\lceil \frac{i}{d} \right\rceil^2 = \frac{\beta \pi^2}{2T^2} d^{-3} N^{-3} \sum_{i=1}^N d i^2 \to \frac{\beta \pi^2}{6T^2} d^{-2},
$$

so we have found $\tau^{\mathrm{RWM}}$. A similar calculation gives the required limit $\tau^{\mathrm{SLA}}$ for SLA. $\square$

Proof of Proposition 2. In this case $\phi_n(x) = \sum_{i=1}^N G(X_t) \Delta t$, for values $X_t = X_t(x)$ from (4.9). We adjust appropriately the specification of the norms. So, for $X = (X_t)_{t=1}^N$ with $X_t \in \mathbb{R}^d$ we define

$$
\|X\|_{L_p} = \left( \sum_{t=1}^N |X_t|^p \Delta t \right)^{1/p}, \qquad \|X\|_s = \left( \sum_{i=1}^N i^{2s} |x_{\cdot,i}|^2 \right)^{1/2}.
$$

Let $x = (x_{\cdot,1}^\top, x_{\cdot,2}^\top, \dots, x_{\cdot,N}^\top)$, $y = (y_{\cdot,1}^\top, y_{\cdot,2}^\top, \dots, y_{\cdot,N}^\top)$ be elements of $\mathbb{R}^n$, for $n = d \times N$, and $X = (X_t)_{t=1}^N$, $(Y_t)_{t=1}^N$ the corresponding discrete-time paths from (4.9). Then, working exactly as in (4.7) and (4.8), using the discrete version of the Sobolev embedding in Lemma B.3(ii), we obtain that

$$
|\phi_n(y) - \phi_n(x)| \leq C(1 + \|X\|_s^p + \|Y - X\|_s^p)\|Y - X\|_{L_2}
$$

for some $s = s(p) < 1/2$. Note now that $\|X\|_{L_2} \equiv |x|$ and $\|X\|_s \leq C|x|_s$ for arbitrary $X$, thus $\{\phi_n\}$ satisfies Condition 4 for $s = s(p) < 1/2$ and $s' = 0$. In this context $\kappa = 1$. Indeed, under $\tilde{\pi}_n$, $x \sim \prod_{i=1}^n \mathcal{N}(0, \Lambda_{i,n}^2)$ where

$$
\Lambda_{i,n} = \sqrt{\frac{1}{2\beta}} T \left( \sin\left( \frac{\lceil i/d \rceil \pi}{2(N+1)} \right)(N+1) \right)^{-1}.
$$

Using the fact that $C_- \leq \sin(v\frac{\pi}{2})/v \leq C_+$ when $v \in (0,1)$ for constants $C_-, C_+ > 0$ we get that $C_- i^{-1} \leq \Lambda_{i,n} \leq C_+ i^{-1}$ for some (other) constants $C_-, C_+ > 0$. It remains to identify $\tau^{\text{RWM}}$ and $\tau^{\text{SLA}}$. For the case of RWM:

$$n^{-2\kappa-1}\sum_{i=1}^{n}\Lambda_{i,n}^{-2} = d^{-3}N^{-3}\frac{2\beta}{T^2}\sum_{i=1}^{N}d\sin^2\left(\frac{i\pi}{2(N+1)}\right)(N+1)^2$$

$$\rightarrow \frac{2\beta}{T^2}d^{-2}\int_0^1 \sin^2\left(v\frac{\pi}{2}\right)dv,$$

so $\tau^{\text{RWM}}$ is as stated in the proposition. One can similarly calculate $\tau^{\text{SLA}}$. $\square$

PROOF OF PROPOSITION 3. For $x, y \in \mathbb{R}^n$ with the structure (5.5), with $n = 4N(N+1)$, let $X_N = P_N(x)$ and $Y_N = P_N(y)$. Recalling the definition of $\sigma = \sigma((i,j))$ we can write for any real $s$,

$$\|X_N\|_s^2 = \sum_{k \in \mathbb{Z}_{U,N}^2} |k|^{2s}((x_k^{\sin})^2 + (x_k^{\cos})^2) = \sum_{i=1}^{n/2}|\sigma^{-1}(i)|^{2s}(x_{2i-1}^2 + x_{2i}^2).$$

From the particular ordering of the elements of $x$ one can easily see that

$$C_- N_0^2 \leq |\sigma^{-1}(i)|^2 \leq C_+ N_0^2, \qquad \text{when } 2(N_0-1)N_0 + 1 \leq i \leq 2N_0(N_0+1)$$

for some $N_0 \leq N$. Now, for a given $i$ and corresponding $N_0 = N_0(i)$ satisfying the second inequality it is true that $C_- i^{1/2} \leq N_0 \leq C_+ i^{1/2}$, therefore

$$(\text{B.12}) \qquad\qquad C_- i \leq |\sigma^{-1}(i)|^2 \leq C_+ i.$$

This gives that

$$C_-\|X_N\|_s \leq |x|_{s/2} \leq C_+\|X_N\|_s.$$

Using now (5.7) we obtain that

$$|\phi_n(y) - \phi_n(x)| = |\phi(Y_N) - \phi(X_N)| \leq \delta(\|X_N\|_{s_0}, \|Y_N\|_{s_0})\|Y_N - X_N\|_{s_0}$$
$$\leq C\delta'(|x|_{s_0/2}, |y|_{s_0/2})|y - x|_{s_0/2}$$

for the locally bounded $\delta' = \delta'(\cdot, \cdot)$ defined as $\delta'(a,b) = \sup_{0 \leq u \leq a, 0 \leq v \leq b} \delta(u,v)$ for $a, b \geq 0$. Again from (5.7),

$$|\phi_n(x)| = |\phi(X_N)| \leq C(1 + \|X_N\|^6) = C(1 + |x|^6).$$

Thus, $\phi_n$ satisfies Condition 5 for parameters $s = s' = s_0/2$ and $s'' = 0$; in this context $\kappa = \alpha/2$ since under $\tilde{\pi}_n$, $x_i \sim \mathcal{N}(0, \Lambda_i^2)$ with

$$\Lambda_i^2 = (4\pi^2)^{-\alpha}|\sigma^{-1}(\lceil i/2 \rceil)|^{-2\alpha}.$$

So, from (B.12), $C_- i^{-\alpha/2} \leq \Lambda_i \leq C_+ i^{-\alpha/2}$. Also,

$$\tau^{\text{RWM}} = \lim_n n^{-(\alpha+1)} \sum_{i=1}^n \Lambda_i^{-2} = 4^{-(\alpha+1)} (4\pi^2)^\alpha \lim_N \left\{ N^{-2(\alpha+1)} 2 \sum_{k \in \mathbb{Z}^2_{U,N}} |k|^{2\alpha} \right\}$$

$$= \frac{1}{2} \pi^{2\alpha} \lim_N \sum_{\substack{-N \leq k_1 \leq N \\ 0 \leq k_2 \leq N}} \left( \frac{k_1^2}{N^2} + \frac{k_2^2}{N^2} \right)^a \frac{1}{N^2},$$

which gives the stated result for $\tau^{\text{RWM}}$. A similar calculation gives $\tau^{\text{SLA}}$. $\quad\square$

**Acknowledgments.** We thank Frank Pinski for providing Figure 1 and Jochen Voss for Figure 2.

## REFERENCES

[1] BÉDARD, M. (2008). Optimal acceptance rates for Metropolis algorithms. *Stochastic Process. Appl.* **118** 2198–2222.

[2] BÉDARD, M. (2007). Weak convergence of Metropolis algorithms for non-i.i.d. target distributions. *Ann. Appl. Probab.* **17** 1222–1244. MR2344305

[3] BÉDARD, M. and ROSENTHAL, J. (2008). Optimal scaling of Metropolis algorithms: Heading toward general target distributions. *Canad. J. Statist.* **36** 483–503.

[4] BENNETT, A. F. (2002). *Inverse Modeling of the Ocean and Atmosphere.* Cambridge Univ. Press, Cambridge. MR1920432

[5] BESKOS, A., ROBERTS, G., STUART, A. and VOSS, J. (2008). MCMC methods for diffusion bridges. *Stoch. Dyn.* **8** 319–350.

[6] BESKOS, S. and STUART, A. (2008). MCMC methods for sampling function space. In *Plenary Lectures Volume, ICIAM 2007.* To appear.

[7] BILLINGSLEY, P. (1999). *Convergence of Probability Measures,* 2nd ed. Wiley, New York. MR1700749

[8] BREYER, L. A., PICCIONI, M. and SCARLATTI, S. (2004). Optimal scaling of MaLa for nonlinear regression. *Ann. Appl. Probab.* **14** 1479–1505. MR2071431

[9] BREYER, L. A. and ROBERTS, G. O. (2000). From Metropolis to diffusions: Gibbs states and optimal scaling. *Stochastic Process. Appl.* **90** 181–206. MR1794535

[10] CHRISTENSEN, O. F., ROBERTS, G. O. and ROSENTHAL, J. S. (2005). Scaling limits for the transient phase of local Metropolis–Hastings algorithms. *J. R. Stat. Soc. Ser. B Stat. Methodol.* **67** 253–268. MR2137324

[11] COTTER, S., DASHTI, M., ROBINSON, J. and STUART, A. (2008). Mathematical foundations of data assimilation problems in fluid mechanics. Unpublished manuscript.

[12] DA PRATO, G. and ZABCZYK, J. (1992). *Stochastic Equations in Infinite Dimensions. Encyclopedia of Mathematics and Its Applications* **44**. Cambridge Univ. Press, Cambridge. MR1207136

[13] DELYON, B. and HU, Y. (2006). Simulation of conditioned diffusion and application to parameter estimation. *Stochastic Process. Appl.* **116** 1660–1675. MR2269221

[14] GEYER, C. J. and THOMPSON, E. A. (1992). Constrained Monte Carlo maximum likelihood for dependent data. *J. R. Stat. Soc. Ser. B Stat. Methodol.* **54** 657–683. With discussion and a reply by the authors. MR1185217

[15] Gilks, W. R., Richardson, S. and Spiegelhalter, D. J., eds. (1996). *Markov Chain Monte Carlo in Practice.* Chapman & Hall, London. MR1397966

[16] Hairer, M., Stuart, A. M., Voss, J. and Wiberg, P. (2005). Analysis of SPDEs arising in path sampling. I. The Gaussian case. *Commun. Math. Sci.* **3** 587–603. MR2188686

[17] Hall, P. and Heyde, C. C. (1980). *Martingale Limit Theory and Its Application.* Academic Press, New York. MR624435

[18] Hastings, W. K. (1970). Monte Carlo sampling methods using Markov chains and their applications. *Biometrica* **57** 97–109.

[19] Kalnay, E. (2003). *Atmospheric Modeling, Data Assimilation and Predictability.* Cambridge Univ. Press, Cambridge.

[20] Neal, P. and Roberts, G. (2006). Optimal scaling for partially updating MCMC algorithms. *Ann. Appl. Probab.* **16** 475–515. MR2244423

[21] Roberts, G. O., Gelman, A. and Gilks, W. R. (1997). Weak convergence and optimal scaling of random walk Metropolis algorithms. *Ann. Appl. Probab.* **7** 110–120. MR1428751

[22] Roberts, G. O. and Rosenthal, J. S. (1998). Optimal scaling of discrete approximations to Langevin diffusions. *J. R. Stat. Soc. Ser. B Stat. Methodol.* **60** 255–268. MR1625691

[23] Robinson, J. C. (2001). *Infinite-Dimensional Dynamical Systems: An Introduction to Dissipative Parabolic PDEs and the Theory of Global Attractors.* Cambridge Univ. Press, Cambridge. MR1881888

[24] Sherlock, C. and Roberts, G. (2008). Optimal scaling of the random walk Metropolis on elliptically symmetric unimodal targets. *Bernoulli.* To appear.

A. Beskos
G. Roberts
Department of Statistics
University of Warwick
Coventry, CV4 7AL
United Kingdom
E-mail: a.beskos@warwick.ac.uk
gareth.o.roberts@warwick.ac.uk

A. Stuart
Mathematics Institute
University of Warwick
Coventry, CV4 7AL
United Kingdom
E-mail: A.M.Stuart@warwick.ac.uk